\documentclass[reqno,11pt]{amsart}

\usepackage{fullpage}
\usepackage{amssymb}
\usepackage{enumitem}
\usepackage[foot]{amsaddr}
\usepackage{placeins}
\usepackage[hidelinks]{hyperref}

\allowdisplaybreaks

\theoremstyle{plain}
\newtheorem{theorem}{Theorem}[section]
\newtheorem{proposition}[theorem]{Proposition}
\newtheorem{corollary}[theorem]{Corollary}
\newtheorem{conjecture}[theorem]{Conjecture}

\theoremstyle{definition}
\newtheorem{definition}[theorem]{Definition}
\newtheorem{question}[theorem]{Question}

\theoremstyle{remark}
\newtheorem{remark}[theorem]{Remark}

\title{Groups of Prime Units}
\author{Trey Li}

\address{The University of Manchester, Manchester, UK}
\email{trey.li@manchester.ac.uk}

\date{October 29, 2026}



\begin{document}

\begin{abstract}
We study groups of units modulo $n$ whose nonidentity elements are all primes in $\{1,\ldots,n-1\}$. We prove that there are only nine groups of prime units which are the full group of units modulo $n$, and there are fewer than $10^{22}$ groups of prime units containing all prime units in the possible range $(\sqrt n,n-\sqrt n)\cup\{n-1\}$. We prove under Dickson's conjecture that there are  infinitely many groups of prime units of arbitrary order $d \ge 2$ modulo composite integers. We prove under Schinzel's Hypothesis H that there are infinitely many groups of prime units of order 3 modulo a prime $r$. We prove that the number of such prime moduli $r$ smaller than $X$ is $\ll X/(\log X)^3$.\end{abstract}

\maketitle

\tableofcontents

\section{Introduction}

We define a group of prime units modulo $n$ to be a group of units modulo $n$ whose nonidentity elements are all primes smaller than $n$. We write $\mathbb P_{n,d}^{\times}$ for such a group of order $d$.

This definition only makes sense for the standard representatives in $\{1,\ldots,n-1\}$, since by Dirichlet's theorem every invertible residue class contains infinitely many primes \cite{Davenport2000} and thus if arbitrary representatives were allowed, then every group of units could be represented by $1$ and primes. 

One might expect such groups to be rare, since requiring a group to be full of primes seems restrictive. However, a simple observation shows that the number of such groups is infinite. In fact, for any prime $p$, the set $\{1,p\}$ is always a group of prime units modulo $p+1$ since $p^2\equiv1\pmod{p+1}$. 

But this example is not surprising, as each group $\{1,p\}$ has only one prime. One might then guess that such groups cannot be very large. This naturally leads us to consider two special types of groups containing as many primes as possible for the same modulus. A group of prime units modulo $n$ is called \textit{complete} if it is the full unit group modulo $n$, and \textit{compact} if it contains every prime unit in $(\sqrt n,n-\sqrt n)\cup\{n-1\}$. This range $(\sqrt n,n-\sqrt n)\cup\{n-1\}$ is the only possible range for the nonidentity elements of such a group, since any prime outside this range would force the group to contain a composite element. We show that, for these two types of groups, the intuition that the group sizes cannot be very large is indeed correct.

\begin{theorem} \label{thm:complete-compact}
The order of a complete group of prime units is at most $8$.
The order of a compact group of prime units is less than $10^{22}$.
\end{theorem}

Having seen that this intuition holds for complete and compact groups, one might expect that even groups that are neither complete nor compact cannot be very large. However, under Dickson's conjecture, we show that they can be arbitrarily large. 

\begin{theorem}\label{thm:general-orders}
If Dickson's conjecture holds, then for every integer $d\geq 2$, there are infinitely many composite moduli
$n$ for which a group of prime units of order $d$ exists.\end{theorem}

The proof of this theorem uses the Chinese remainder theorem to construct groups modulo composite integers whose nonidentity elements are given by linear polynomials, and Dickson's conjecture
\cite{Dickson1904} to make them all prime.

The theorem above concerns composite moduli. For prime moduli, the problem is more restrictive. To see this, notice that every group of prime units modulo $n$ is a subgroup of the group of units modulo $n$. For a prime modulus $r$, since the group of units modulo $r$ is cyclic, it has exactly one subgroup of order $d$ for each $d\mid r-1$. Thus, a group of prime units of order $d$ modulo $r$ is unique if it exists. Moreover, for $r>3$, no subgroup of even order can be a group of prime units, since it contains $-1$, whose standard representative $r-1$ is composite. Hence, proving the existence of infinitely many groups of prime units for prime moduli is less straightforward than for composite moduli. 

We study the problem order by order. Since groups of order $1$ are trivial and order $2$ occurs only modulo $3$, the first nontrivial case is order $3$. We focus on this case and call such groups \emph{trey groups} and their prime moduli \emph{trey primes}.

The order-$3$ case has a simple connection to quadratic congruences. For a prime $r>3$, the two nonidentity elements of a subgroup of order $3$ are exactly the two roots of $x^2+x+1\equiv0\pmod r$. Hence $r$ is a trey prime if and only if both roots exist and are distinct primes. Roots of quadratic congruences modulo primes have been considered by Duke, Friedlander, and Iwaniec \cite{DFI95}. Our study of trey groups asks whether both standard representatives of these roots are prime and seeks quantitative bounds for the corresponding prime moduli $r$.

We show that the number of trey primes is infinite if Schinzel's Hypothesis H holds.

\begin{theorem}\label{thm:schinzel}
If Schinzel's Hypothesis H holds, then there are infinitely many trey primes.
\end{theorem}

A further question is how sparse trey primes are. Let $N_3(X)$ be the number of trey primes up to $X$. More generally, for fixed $d$, let $N_d(X)$ count the prime moduli $r\leq X$ such that $\mathbb P_{r,d}^{\times}$ exists.

\begin{theorem}\label{thm:main-upper}
For all $X\geq2$, we have $N_3(X)\ll\frac{X}{(\log X)^3}$.
\end{theorem}

This is the main theorem of the paper. The proof of this theorem uses the relations in a trey group to construct a binary quadratic form, which allows us to represent each trey group by an integer matrix. We then apply the affine sieve of Bourgain, Gamburd, and Sarnak \cite{BGS10,BGS11} to count matrices satisfying the three prime conditions. By computing the local densities, we show that the sieve has dimension $3$. This connects our problem to the study of prime values on group orbits \cite{HoreshNevo,NevoSarnak}.

The rest of the paper is organized as follows.
Section~2 studies complete and compact groups of prime units and proves Theorem~\ref{thm:complete-compact}.
Section~3 proves Theorem~\ref{thm:general-orders}.
Section~4 studies trey groups and trey primes and proves Theorem~\ref{thm:schinzel}.
Section~5 proves Theorem~\ref{thm:main-upper} using binary quadratic forms and the affine sieve.
Section~6 presents open problems.

\section{Finitely many complete and compact groups of prime units}\label{sec:elementary} 

We define groups of prime units as follows.

\begin{definition}
Let $n\geq2$ and $d\geq1$. A \emph{group of prime units} of order $d$ modulo $n$, denoted by $\mathbb P_{n,d}^{\times}$, is a subgroup of $(\mathbb Z/n\mathbb Z)^\times$ of order $d$ whose nonidentity elements have prime standard representatives in $\{1,\ldots,n-1\}$.
\end{definition}

For a prime modulus $n = r$, the group $\mathbb P_{r,d}^{\times}$ is unique if it exists, since $(\mathbb Z/r\mathbb Z)^\times$ is cyclic and there is a unique subgroup of order $d$ for each $d \mid r-1$. However, for a composite modulus $n$, the notation $\mathbb P_{n,d}^{\times}$ can mean more than one group. For example, $\mathbb P_{8,2}^{\times}$ can mean $\{1,3\}$, $\{1,5\}$, or $\{1,7\}$.

For moduli $n = 2,4,p^{j}$ or $2p^{j}$, where $p$ is an odd prime and $j\geq1$, the group $\mathbb P_{n,d}^\times$ is cyclic. This is because every group $(\mathbb Z/n\mathbb Z)^\times$ with $n =$ $2$, $4$, $p^j$, or $2p^j$ is cyclic, every subgroup of a cyclic group is cyclic, and every group $\mathbb P_{n, d}^\times$ is a subgroup of $(\mathbb Z/n\mathbb Z)^\times$.

For moduli $n \ne 2,4,p^j,2p^j$, two groups of prime units are not guaranteed to be isomorphic even when they have the same modulus $n$ and order $d$. For example, when $n = 30$ and $d = 4$, the group $\{1,7,13,19\}$ modulo $30$ is cyclic but the group $\{1,11,19,29\}$ modulo $30$ is a Klein four-group. Even if the two groups of prime units share the same integers as their elements, they can still be non-isomorphic if their moduli are different. For example, the same set $\{1,7,13,19\}$ forms a cyclic group modulo $30$ but a Klein four-group modulo $24$.

The nonidentity elements of a group of prime units lie either in the middle range of $\{1,\ldots,n-1\}$ or at $n-1$.

\begin{proposition}\label{prop:general-range}
For $n\geq2$, every nonidentity element of a group of prime units modulo $n$ lies in $(\sqrt n,n-\sqrt n) \cup \{n-1\}$.
\end{proposition}

\begin{proof}
Suppose $p\leq\sqrt n$. We have $p^2 \le n$. But $p$ is coprime to $n$. Thus $p^2\ne n$ and $p^2< n$. So the group element $p^2 \mod n = p^2$ is composite, a contradiction. Hence $p>\sqrt n$.

Now suppose $n-\sqrt n\leq p<n-1$. Write $p=n-a$, where $1<a\leq\sqrt n$. We have $p^2\equiv a^2\pmod n$ and $1<a^2\le n$. But $a$ is coprime to $n$. Thus $a^2\ne n$ and $1<a^2< n$. So the group element $p^2 \mod n = a^2$
is composite, a contradiction. Hence $p\in(\sqrt n,n-\sqrt n)\cup\{n-1\}$. \end{proof}

This restriction leads to two special types of groups of prime units.

\begin{definition}
A group of prime units $G$ modulo $n$ is \emph{complete} if $G=(\mathbb Z/n\mathbb Z)^\times$. It is
\emph{compact} if it contains all prime units modulo $n$ in $(\sqrt n,n-\sqrt n)$ and also contains
$n-1$ when $n-1$ is prime.
\end{definition}

Every complete group is compact. This is because a complete group contains all units and thus it must contain whatever a compact group requires to contain.

A compact group may not be complete. For example, $\{1\}$ is compact modulo $5$, since there are no prime units in $(\sqrt 5,5-\sqrt 5)$ and $4$ is not prime. But it is not complete, since $(\mathbb Z/5\mathbb Z)^\times=\{1,2,3,4\}$.

However, whenever both the complete and compact groups exist for the same modulus $n$, they must be the same group. To see this, let $G$ be the complete group and $H$ be the compact group modulo $n$. On the one hand, since $G$ is complete and all elements of $H$ are units, we have $H\subseteq G$. On the other hand, since $H$ is compact and all nonidentity elements of $G$ are in $(\sqrt n,n-\sqrt n)$ or equal $n-1$ by Proposition~\ref{prop:general-range}, we have $G\subseteq H$. So $G=H$.

There are only finitely many complete and compact groups. In fact, there are exactly 9 complete groups and at most $10^{22}$ compact groups.

\begin{theorem}\label{thm:complete-classification}
There are only finitely many complete groups of prime units.
\end{theorem}

\begin{proof}
We prove that no complete group of prime units exists modulo $n>120$. Suppose for a contradiction that such a group exists for $n > 120$. Let $\mathbb P_n^\times=(\mathbb Z/n\mathbb Z)^\times=\{1,q_1,\ldots,q_k\}$, where $q_1$ is the smallest nonidentity element.

Since $q_1$ is in the group, so is $q_1^2$ modulo $n$. We must have $n\le q_1^2$ otherwise the group would contain the composite number $q_1^2$ as an element. Hence for $n>120$, we must have $q_1\geq11$.

Let $p_1,\ldots,p_s$ be all the primes smaller than $q_1$. Namely $q_1 = p_{s+1}$. Since $q_1\geq11$, we have $s\geq4$. Also, every $p_i$ must divide $n$, because otherwise it would be a unit modulo $n$ and hence should be in the group as well, contradicting the choice of $q_1$. Since these primes are different, their product also divides $n$. We therefore have $p_1\cdots p_s\leq n\leq q_1^2=p_{s+1}^2$.

To obtain a contradiction, we show that the opposite also holds, namely $p_1\cdots p_s>p_{s+1}^2$ for all $s\geq4$. We prove this by induction.

For $s=4$, the inequality $p_1\cdots p_s>p_{s+1}^2$ holds because $2\cdot3\cdot5\cdot7=210>11^2$. Now suppose that $p_1\cdots p_s>p_{s+1}^2$ for some $s\geq4$. We want to prove that $p_1\cdots p_{s+1}>p_{s+2}^2$.

Multiplying the induction hypothesis $p_1\cdots p_s>p_{s+1}^2$ by $p_{s+1}$ gives $p_1\cdots p_{s+1}>p_{s+1}^3$. Since $p_{s+1}>4$, we have $p_{s+1}^3>4p_{s+1}^2$. Moreover, Bertrand's postulate gives $p_{s+2}<2p_{s+1}$, so $4p_{s+1}^2>p_{s+2}^2$. Combining these inequalities, we have $p_1\cdots p_{s+1}>p_{s+2}^2$, as required. This completes the induction.

Therefore, if there exists a complete group of prime units modulo $n>120$, then for every $s\geq4$, we have both $p_1\cdots p_s\leq p_{s+1}^2$ and $p_1\cdots p_s>p_{s+1}^2$, a contradiction. Hence no complete group of prime units exists modulo $n>120$.

Since there are only finitely many integers $n\leq120$, and each modulus has at most one complete group, the theorem follows.
\end{proof}

A direct check of all $n\leq120$ gives all possible complete groups of prime units as follows.
$$
\begin{array}{c|l}
n&(\mathbb Z/n\mathbb Z)^\times\\ \hline
2&\{1\}\\
3&\{1,2\}\\
4&\{1,3\}\\
6&\{1,5\}\\
8&\{1,3,5,7\}\\
12&\{1,5,7,11\}\\
18&\{1,5,7,11,13,17\}\\
24&\{1,5,7,11,13,17,19,23\}\\
30&\{1,7,11,13,17,19,23,29\}.
\end{array}
$$

\begin{theorem}\label{thm:compact-finiteness}
There are only finitely many compact groups of prime units.
\end{theorem}

\begin{proof}
We prove that no compact group exists modulo $n\geq10^{22}$. Since each modulus has at most one compact group, this will prove the theorem.

Suppose for contradiction that a compact group exists modulo $n\geq10^{22}$. We will show that this group must contain a composite number, which contradicts the definition of groups of prime units. For this, it suffices to show that there exist two prime units $p$ and $q$ in the interval $(\sqrt n,n-\sqrt n)$ hence in the compact group such that $pq$ modulo $n$ is composite. Notice that if $n<pq<2n$, then $pq$ modulo $n$ is $pq-n$. Hence it suffices to find $p$ and $q$ such that $n<pq<2n$ and $pq-n$ is composite. To ensure $n<pq<2n$, we can choose $p$ and $q$ from the interval $(11\sqrt n/10,6\sqrt n/5)\subset(\sqrt n,n-\sqrt n)$, which gives $n<121n/100<pq<36n/25<2n$. To ensure that $pq-n$ is composite, it suffices to find a prime $\ell$ such that $p\equiv1\pmod\ell$, $q\equiv n\pmod\ell$, and $pq-n>\ell$. The first two conditions imply $pq\equiv n\pmod\ell$ and thus $\ell\mid(pq-n)$, while the third ensures that $pq - n \ne \ell$ and thus $\ell$ is a proper divisor of $pq-n$ and thus $pq-n$ is composite. We can ensure $pq-n>\ell$ by requiring $\ell\leq2\log n$, since $pq-n>21n/100>2\log n$ for $n\geq10^{22}$. Hence our goal reduces to finding prime units $p,q\in(11\sqrt n/10,6\sqrt n/5)$ and a prime $\ell\leq2\log n$ such that $p\equiv1\pmod\ell$ and $q\equiv n\pmod\ell$.

To find the required primes $p$ and $q$, we first need to choose $\ell$, since the congruences $p\equiv1\pmod\ell$ and $q\equiv n\pmod\ell$ depend on $\ell$. We will use a theorem that counts primes in residue classes modulo $\ell$, but this theorem requires the residue classes to be coprime to $\ell$. In our case, the residue classes are $1$ and $n$. Since $1$ is always coprime to $\ell$, we only need to ensure that $n$ is coprime to $\ell$, which means $\ell\nmid n$ since $\ell$ is prime. Therefore, we first prove that there exists an odd prime $\ell\leq2\log n$ not dividing $n$. Once such a prime $\ell$ is fixed, we can use the theorem to find the required primes $p$ and $q$.

To prove that an odd prime $\ell\leq2\log n$ not dividing $n$ exists, it suffices to rule out the possibility that every odd prime not exceeding $2\log n$ divides $n$. This is because otherwise if all these primes divided $n$, then their product would also divide $n$ since they are distinct primes, and hence the product could not exceed $n$. Therefore, if we can prove that their product is greater than $n$, then we are sure that at least one of these primes cannot divide $n$. We can then choose that prime as $\ell$.

To prove that the product of all odd primes not exceeding $2\log n$ is greater than $n$, we take the logarithm of this product, which equals the sum of $\log r$ over all odd primes $r\leq2\log n$. We use the known bound~\cite{Eberl2024} that the sum of $\log r$ over all primes $r\leq t$ is at least $0.82t$ whenever $t\geq97$, namely $\vartheta(t)\geq0.82t$ for $t\geq97$, where $\vartheta(t)=\sum_{p\leq t}\log p$ is the first Chebyshev function \cite{Eberl2024}. Since $n\geq10^{22}$, we have $2\log n>97$, so the sum of $\log r$ over all primes $r\leq2\log n$ is at least $0.82(2\log n)=1.64\log n$. To obtain the sum over odd primes only, we subtract $\log2$, which is the contribution of the prime $2$. Hence the logarithm of the product of all odd primes not exceeding $2\log n$ is at least $1.64\log n-\log2>\log n$. Therefore, the product of these odd primes is greater than $n$. Since their product exceeds $n$, not all these primes can divide $n$. We can thus choose an odd prime $\ell\leq2\log n$ that does not divide $n$.

We now fix this prime $\ell$ and turn to finding $p$ and $q$. To find prime units $p,q\in(11\sqrt n/10,6\sqrt n/5)$ such that $p\equiv1\pmod\ell$ and $q\equiv n\pmod\ell$, it suffices to show that this interval contains at least two primes congruent to $1$ modulo $\ell$ and at least two primes congruent to $n$ modulo $\ell$. We require at least two primes in each residue class because every prime in this interval is greater than $\sqrt n$, so at most one prime in the interval can divide $n$ since two distinct prime divisors greater than $\sqrt n$ would have a product greater than $n$. If such a prime divisor of $n$ exists, it is not a unit and hence we cannot choose it as $p$ or $q$. So we need to exclude it. After excluding this prime, if it exists, at least one prime unit remains in each residue class. We can then choose $p$ and $q$ from the two residue classes, respectively. Note that the two residue classes may coincide, and the selected primes $p$ and $q$ need not be distinct, but these do not affect the argument.

To guarantee at least two primes congruent to $1$ modulo $\ell$ and at least two primes congruent to $n$ modulo $\ell$ in the open interval $(11\sqrt n/10,6\sqrt n/5)$, it suffices to show that the half-open interval $(11\sqrt n/10,6\sqrt n/5]$ contains at least three primes satisfying each congruence. We require three primes in each residue class because $6\sqrt n/5$ might itself be prime and we need to exclude it for counting primes in the open interval $(11\sqrt n/10,6\sqrt n/5)$. Excluding this number removes at most one prime from each residue class, leaving at least two primes in each residue class, as required.

To prove that $(11\sqrt n/10,6\sqrt n/5]$ contains at least three primes satisfying each of the congruences $r\equiv1\pmod\ell$ and $r\equiv n\pmod\ell$, it suffices to show that the number of primes satisfying each congruence in this interval is greater than $2$. For either congruence, this number equals the number of primes satisfying it up to $6\sqrt n/5$, minus the number satisfying it up to $11\sqrt n/10$. Therefore, it suffices to find a lower bound for the first count and an upper bound for the second count, and then prove that the difference between these bounds is greater than $2$.

To obtain these bounds, we use Theorem 1.3 of Bennett, Martin, O'Bryant, and Rechnitzer~\cite{BennettMartinOBryantRechnitzer}. The theorem states that, for any integer $a$ coprime to $\ell$, the number of primes $r\leq y$ satisfying $r\equiv a\pmod\ell$ differs from $\frac{1}{\ell-1}\int_2^y\frac{dt}{\log t}$ by less than $\frac{y}{160(\log y)^2}$, provided that $y\geq8\cdot10^9$ when $3\leq\ell\leq10^5$, or $y\geq\exp(0.03\sqrt\ell(\log\ell)^3)$ when $\ell>10^5$. Since $\ell\nmid n$, both $1$ and $n$ are coprime to $\ell$. Thus, once we verify that $11\sqrt n/10$ and $6\sqrt n/5$ satisfy the theorem's condition on $y$, we can apply the theorem to both congruences at both endpoints.

To use Theorem 1.3 with $a=1$ and $a=n$, we must first verify that the theorem applies. Since $\ell\nmid n$, both $1$ and $n$ are coprime to $\ell$. We therefore only need to verify that $11\sqrt n/10$ and $6\sqrt n/5$ satisfy the lower bound on $y$ required by the theorem. Since both numbers are greater than $\sqrt n$, it suffices to verify the lower bound for $\sqrt n$.

If $\ell\leq10^5$, then $\sqrt n\geq10^{11}>8\cdot10^9$, so the required lower bound is satisfied.

If $\ell>10^5$, the theorem requires $y\geq\exp(0.03\sqrt\ell(\log\ell)^3)$. Since both $11\sqrt n/10$ and $6\sqrt n/5$ are greater than $\sqrt n$, it suffices to prove that $\sqrt n>\exp(0.03\sqrt\ell(\log\ell)^3)$, or equivalently, $\log\sqrt n>0.03\sqrt\ell(\log\ell)^3$. Since $\ell\leq4\log\sqrt n$, we have $\sqrt\ell\leq\sqrt{4\log\sqrt n}$ and $\log\ell\leq\log(4\log\sqrt n)$. Therefore,
$$
0.03\sqrt\ell(\log\ell)^3\leq0.03\sqrt{4\log\sqrt n}(\log(4\log\sqrt n))^3.
$$
Hence, to prove that $\log\sqrt n>0.03\sqrt\ell(\log\ell)^3$, it suffices to show that
$$
\log\sqrt n>0.03\sqrt{4\log\sqrt n}(\log(4\log\sqrt n))^3.
$$
Now $\ell>10^5$ and $\ell\leq4\log\sqrt n$ imply $\log\sqrt n>25000$. For $\log\sqrt n\geq25000$, the ratio of the left-hand side to the right-hand side is
$$
\frac{\sqrt{\log\sqrt n}}{0.06(\log(4\log\sqrt n))^3}.
$$
This ratio is increasing because $\log(4\log\sqrt n)>6$, which is because $\log\sqrt n>25000$ implies $4\log\sqrt n>10^5>e^6$; and the value of this ratio at $\log\sqrt n=25000$ is greater than $1.7>1$. Therefore, the required inequality holds, and both endpoints $11\sqrt n/10$ and $6\sqrt n/5$ satisfy the theorem's lower bound on $y$ when $\ell>10^5$.

We can now apply Theorem 1.3 of Bennett, Martin, O'Bryant, and Rechnitzer~\cite{BennettMartinOBryantRechnitzer}. For $a=1$ or $a=n$, the number of primes $r$ satisfying $11\sqrt n/10<r\leq6\sqrt n/5$ and $r\equiv a\pmod\ell$ is greater than
$$
\begin{aligned}
&\left(\frac{1}{\ell-1}\int_2^{6\sqrt n/5}\frac{dt}{\log t}
-\frac{6\sqrt n}{800(\log(6\sqrt n/5))^2}\right) -\left(\frac{1}{\ell-1}\int_2^{11\sqrt n/10}\frac{dt}{\log t}
+\frac{11\sqrt n}{1600(\log(11\sqrt n/10))^2}\right).
\end{aligned}
$$
Expanding the brackets, this expression equals
$$
\begin{aligned}
&\frac{1}{\ell-1}\left(\int_2^{6\sqrt n/5}\frac{dt}{\log t}
-\int_2^{11\sqrt n/10}\frac{dt}{\log t}\right) -\frac{6\sqrt n}{800(\log(6\sqrt n/5))^2}
-\frac{11\sqrt n}{1600(\log(11\sqrt n/10))^2}.
\end{aligned}
$$
Since the difference of the two integrals equals $\int_{11\sqrt n/10}^{6\sqrt n/5}\frac{dt}{\log t}$, we obtain the lower bound
$$
\begin{aligned}
&\frac{1}{\ell-1}\int_{11\sqrt n/10}^{6\sqrt n/5}\frac{dt}{\log t} -\frac{6\sqrt n}{800(\log(6\sqrt n/5))^2}
-\frac{11\sqrt n}{1600(\log(11\sqrt n/10))^2}.
\end{aligned}
$$

To prove that there are at least three primes satisfying each congruence, it suffices to show that this lower bound exceeds $2$, since the number of primes is an integer.

Since $\log t\leq\log(6\sqrt n/5)$ for $11\sqrt n/10<t\leq6\sqrt n/5$, we have $\frac{1}{\log t}\geq\frac{1}{\log(6\sqrt n/5)}$. The length of this interval is $6\sqrt n/5-11\sqrt n/10=\sqrt n/10$. Therefore, $\int_{11\sqrt n/10}^{6\sqrt n/5}\frac{dt}{\log t}\geq\frac{\sqrt n}{10\log(6\sqrt n/5)}$. Also, since both $11\sqrt n/10$ and $6\sqrt n/5$ are greater than $\sqrt n$, both logarithms $\log(6\sqrt n/5)$ and $\log(11\sqrt n/10)$ in the two error terms $\frac{6\sqrt n}{800(\log(6\sqrt n/5))^2}$ and $\frac{11\sqrt n}{1600(\log(11\sqrt n/10))^2}$ are greater than $\log\sqrt n$. Hence the sum of the two error terms is less than $\frac{6\sqrt n}{800(\log\sqrt n)^2}+\frac{11\sqrt n}{1600(\log\sqrt n)^2}=\frac{23\sqrt n}{1600(\log\sqrt n)^2}$. Substituting these estimates into the lower bound earlier, the number of primes satisfying $r\equiv a\pmod\ell$ in $(11\sqrt n/10,6\sqrt n/5]$ is greater than
$$
\frac{\sqrt n}{10(\ell-1)\log(6\sqrt n/5)}-\frac{23\sqrt n}{1600(\log\sqrt n)^2}.
$$

To see that this lower bound is greater than $2$, note that the denominator $10(\ell-1)\log(6\sqrt n/5)$ is less than $40.4(\log\sqrt n)^2$, since $\ell-1<4\log\sqrt n$ by $\ell\leq2\log n=4\log\sqrt n$, and $\log(6\sqrt n/5)=\log\sqrt n+\log(6/5)<1.01\log\sqrt n$ because $\log(6/5)<0.01\cdot11\log10\leq0.01\log\sqrt n$ for $n\geq10^{22}$. Therefore, the lower bound is greater than
$$
\left(\frac{1}{40.4}-\frac{23}{1600}\right)\frac{\sqrt n}{(\log\sqrt n)^2}>\frac{\sqrt n}{100(\log\sqrt n)^2},
$$ since $\frac{1}{40.4}-\frac{23}{1600}>\frac{1}{100}$. Now since $\sqrt n\geq10^{11}$ and the function $x/(\log x)^2$ is increasing for $x>e^2$, we have that
$$
\frac{\sqrt n}{100(\log\sqrt n)^2}\geq\frac{10^{11}}{100(11\log10)^2}>2.
$$

Therefore, $(11\sqrt n/10,6\sqrt n/5]$ contains at least three primes congruent to $1$ modulo $\ell$ and at least three primes congruent to $n$ modulo $\ell$. Removing $6\sqrt n/5$, if it is prime, leaves at least two primes satisfying each congruence in $(11\sqrt n/10,6\sqrt n/5)$.

Since at most one prime in this interval divides $n$, we can exclude that prime, if it exists, and still choose prime units $p$ and $q$ satisfying $p\equiv1\pmod\ell$ and $q\equiv n\pmod\ell$.

We have now found the primes required at the beginning of the proof. Since $p$ and $q$ are prime units in $(\sqrt n,n-\sqrt n)$, both belong to the compact group. Their product modulo $n$ must therefore also belong to the group.

However, the representative of $pq$ modulo $n$ is $pq-n$. The congruences $p\equiv1\pmod\ell$ and $q\equiv n\pmod\ell$ imply $\ell\mid(pq-n)$. Moreover, $pq-n>21n/100>2\log n\geq\ell$. Thus $\ell$ is a prime divisor of $pq-n$ smaller than $pq-n$ (i.e., $\ell \ne pq-n$), so $pq-n$ is composite. This contradicts the definition of a group of prime units.

Therefore no compact group exists modulo $n\geq10^{22}$, and the theorem follows.
\end{proof}

From the theorem, we also learn that the order $d$ of a compact group is less than $10^{22}$, since $d\leq\varphi(n)<n<10^{22}$. A more precise estimate can be obtained by counting the primes in the interval $(\sqrt n,n-\sqrt n)$. If $k$ denotes the number of primes in this interval, then $k\leq d\leq k+2$, since at most one of these primes can divide $n$, the identity $1$ always contributes one additional element, and $n-1$ contributes another element if it is prime.

\section{Infinitely many groups of prime units}\label{sec:general}

Having shown that complete and compact groups of prime units are finite in number, we now turn to general groups of prime units. Unlike complete and compact groups, there are infinitely many general groups of prime units. As is mentioned in the Introduction, there are infinitely many groups $\{1,p\}$ of prime units of order $d = 2$. Now we prove an infinity result for orders $d\geq2$ under Dickson's conjecture.

\begin{conjecture}[Dickson's conjecture \cite{Dickson1904}]\label{conj:dickson}
Let $f_1,\ldots,f_k\in\mathbb Z[x]$ be linear polynomials with positive leading coefficients. Suppose no prime divides $f_1(x)\cdots f_k(x)$ for every integer $x$. Then all the polynomials take prime values simultaneously for infinitely many positive integers $x$.
\end{conjecture}

\begin{proof}[Proof of Theorem~\ref{thm:general-orders}]
Fix an order $d\geq2$. We want to prove that infinitely many different composite moduli have a group of prime units of order $d$. For this, it suffices to use a positive integer $x$ to construct a composite modulus $n_x$ and a group $G_x$ of order $d$ modulo $n_x$, such that different values of $x$ give different moduli $n_x$. It then suffices to prove that, for infinitely many values of $x$, all the nonidentity elements of $G_x$ are primes. To prove this, we will express these elements as linear polynomials in $x$ and use Dickson's conjecture to show that they are all prime for infinitely many values of $x$.

Now we construct the groups $G_x$ of order $d$ modulo $n_x$. To do this, we need to find exactly $d$ different residues modulo $n_x$ that form a group under multiplication. Instead of finding these residues directly, we can start with a group $C$ that already has order $d$ and use it to construct $G_x$. To find such a group $C$, it suffices to choose a prime $\ell$ such that $d$ divides $\ell-1$. This is because $\mathbb F_\ell^\times$ is a cyclic group of order $\ell-1$, and a cyclic group of order $\ell-1$ contains a subgroup of order $d$ whenever $d$ divides $\ell-1$. By Dirichlet's theorem~\cite{Davenport2000}, there exists a prime $\ell\equiv1\pmod d$, so such a subgroup exists. We fix this prime $\ell$ and let $C$ be the subgroup of order $d$ in $\mathbb F_\ell^\times$. Represent its elements by integers in $\{1,\ldots,\ell-1\}$.

We have now found a group $C$ of order $d$ modulo the prime $\ell$. We want to use $C$ to construct a group $G_x$ of order $d$ modulo a composite number $n_x$. For every positive integer $x$, we choose
$$
t_x=\ell(2x+1)+1,\qquad n_x=\ell t_x,\qquad
P_h(x)=1+(h-1)t_x\quad(h\in C).
$$
We define $G_x$ to be the set of residues modulo $n_x$ represented by the numbers $P_h(x)$ for all $h\in C$. Since $t_x>1$, the number $n_x=\ell t_x$ is composite. Also, $n_x$ is strictly increasing with $x$, so different positive integers $x$ give different composite moduli.

To prove that $G_x$ is a group of order $d$, it suffices to show that $G_x$ is isomorphic to $C$. For this, it suffices to show that the correspondence $h\mapsto P_h(x)\pmod{n_x}$ preserves multiplication and that different elements $h\in C$ give different residues $P_h(x)$ modulo $n_x$. These two conditions ensure that the correspondence is an isomorphism from $C$ to $G_x$.

To prove that the correspondence from $C$ to $G_x$ preserves multiplication, it suffices to show that $P_h(x)P_k(x)\equiv P_{hk}(x)\pmod{n_x}$ for every $h,k\in C$, where $hk$ is the product of $h$ and $k$ in $C$. Since $n_x=\ell t_x$, by the Chinese remainder theorem, it suffices to prove that $\ell$ and $t_x$ are coprime as well as $P_h(x)P_k(x)\equiv P_{hk}(x)\pmod{\ell}$ and $P_h(x)P_k(x)\equiv P_{hk}(x)\pmod{t_x}$.

To prove these two congruences, note that $t_x=\ell(2x+1)+1$ gives $t_x\equiv1\pmod\ell$, so $(\ell,t_x)=1$. Since $P_h(x)=1+(h-1)t_x$, we have
$$
P_h(x)\equiv h\pmod\ell,\qquad P_h(x)\equiv1\pmod{t_x}.
$$
Therefore, for any $h,k\in C$, we have $P_h(x)P_k(x)\equiv hk\equiv P_{hk}(x)\pmod\ell$. Similarly, since $P_h(x)$, $P_k(x)$, and $P_{hk}(x)$ are all congruent to $1$ modulo $t_x$, we have $P_h(x)P_k(x)\equiv P_{hk}(x)\pmod{t_x}$. Since $(\ell,t_x)=1$, the Chinese remainder theorem gives $P_h(x)P_k(x)\equiv P_{hk}(x)\pmod{n_x}$. Thus the correspondence preserves multiplication.

To prove that different elements of $C$ give different residues modulo $n_x$, suppose $h,k\in C$ are different. Since the elements of $C$ are represented by different integers in $\{1,\ldots,\ell-1\}$, we have $h\not\equiv k\pmod\ell$. As $P_h(x)\equiv h\pmod\ell$ and $P_k(x)\equiv k\pmod\ell$, the numbers $P_h(x)$ and $P_k(x)$ cannot be congruent modulo $\ell$, and hence cannot be congruent modulo $n_x$.

We have now proved that the correspondence preserves multiplication and sends different elements of $C$ to different elements of $G_x$. Since $G_x$ consists of exactly the residues $P_h(x)$ corresponding to the elements $h$ of $C$, this correspondence is an isomorphism. Therefore, $G_x$ is a group of order $d$.

We have now constructed a group $G_x$ of order $d$ modulo $n_x$. To prove that $G_x$ is a group of prime units for infinitely many positive integers $x$, it remains to show that the standard representatives of all its nonidentity elements are prime for infinitely many $x$. Since $P_1(x)=1$, it represents the identity element. For every $h\neq1$, we have $2\leq h\leq\ell-1$, because the elements of $C$ are represented by integers in $\{1,\ldots,\ell-1\}$. Therefore, since $t_x>1$, we have $1<1+t_x\leq P_h(x)=1+(h-1)t_x\leq1+(\ell-2)t_x<\ell t_x=n_x$. Thus every $P_h(x)$ with $h\neq1$ lies in $\{2,\ldots,n_x-1\}$ and is the standard representative of a nonidentity element of $G_x$ modulo $n_x$. Therefore, it remains to prove that all $P_h(x)$ with $h\neq1$ are prime for infinitely many positive integers $x$.

To prove that the polynomials $\{P_h(x)\}_{h\in C}$ are simultaneously prime for infinitely many $x$, we use Dickson's conjecture. To apply it, it suffices to show that every $P_h(x)$ with $h\neq1$ is a linear polynomial with positive leading coefficient, and that for every prime $s$, there exists an integer $x$ such that none of these polynomials is divisible by $s$.

To verify the first condition, we write
$$
P_h(x)=2\ell(h-1)x+\ell(h-1)+h.
$$
Since $h\neq1$ and $h\in\{1,\ldots,\ell-1\}$, we have $h>1$, so the coefficient $2\ell(h-1)$ of $x$ is positive. Thus every $P_h(x)$ with $h\neq1$ is a linear polynomial with positive leading coefficient.

To verify the second condition, we must show that for every prime $s$, we can choose an integer $x$ such that none of the polynomials $P_h(x)$ with $h\neq1$ is divisible by $s$. We consider three cases.

If $s=2$, none of these polynomials $P_h(x)=1+(h-1)t_x$ is divisible by $s$, since $t_x=\ell(2x+1)+1$ is even with $\ell$ odd. Hence every $P_h(x)$ is odd.

If $s=\ell$, none of these polynomials $P_h(x)=1+(h-1)t_x$ is divisible by $s$, since $P_h(x)\equiv h\pmod\ell$ and every $h\in C$ is nonzero modulo $\ell$.

If $s$ is neither $2$ nor $\ell$, it suffices to choose an integer $x$ such that $t_x\equiv0\pmod s$. This is because $P_h(x)=1+(h-1)t_x$ would then satisfy $P_h(x)\equiv1\pmod s$ for every $h\in C$, so none of these polynomials would be divisible by $s$. Such an integer $x$ exists because $t_x = \ell(2x+1)+1 = 2\ell x + \ell + 1 \equiv0\pmod s$ has an integer solution. To see this, transform the equation into $2\ell x\equiv-(\ell+1)\pmod s$, and notice that a linear congruence $ax\equiv b\pmod s$ has an integer solution whenever $(a,s)=1$. Now $(2\ell,s) = 1$ since $s$ is a prime different from $2$ and $\ell$. Therefore, $2\ell x\equiv-(\ell+1)\pmod s$ must have an integer solution $x$.

Therefore, for every prime $s$, there exists an integer $x$ such that none of the polynomials $P_h(x)$ with $h\neq1$ is divisible by $s$. Hence both conditions of Dickson's conjecture are satisfied. By Dickson's conjecture, there are infinitely many positive integers $x$ for which all $P_h(x)$ with $h\neq1$ are prime. Also, for every such $x$, all nonidentity elements $P_h(x)$ of $G_x$ are primes in $\{2,\ldots,n_x-1\}$, so $G_x$ is a group of prime units of order $d$ modulo $n_x$. Since the numbers $n_x=\ell(\ell(2x+1)+1)$ are composite and strictly increasing with $x$, these groups occur modulo infinitely many different composite moduli $n_x$. This proves the theorem.
\end{proof}

For example, we can construct infinitely many groups of units of order $d=3$ as follows. We choose $\ell=7$ and $C=\{1,2,4\}$, which is a subgroup of order $3$ in $\mathbb F_7^\times$. For every positive integer $x$, our construction gives $t_x=14x+8$ and $n_x=7t_x=98x+56$. The corresponding group is $G_x=\{1,14x+9,42x+25\}$ modulo $n_x$, since $P_1(x)=1$, $P_2(x)=14x+9$, and $P_4(x)=42x+25$. As proved above, $G_x$ is a group of order $3$ for every positive integer $x$. Moreover, Dickson's conjecture implies that there are infinitely many positive integers $x$ for which both $14x+9$ and $42x+25$ are prime. For every such $x$, the group $G_x$ is a group of prime units of order $3$ modulo the composite number $n_x=98x+56$. Since these moduli are different for different values of $x$, this gives infinitely many composite moduli having groups of prime units of order $3$, assuming Dickson's conjecture. For instance, taking $x=1$ gives $n_1=154$ and $G_1=\{1,23,67\}$, where both $23$ and $67$ are primes in $\{2,\ldots,153\}$.

\section{Infinitely many trey groups and trey primes}\label{sec:trey}

We now turn from composite moduli to prime moduli. For a prime $r$, the group $(\mathbb Z/r\mathbb Z)^\times$ is cyclic of order $r-1$, so it has exactly one subgroup of order $d$ for each $d\mid r-1$. Therefore, there is at most one group of prime units of order $d$ modulo $r$. Moreover, when $r=3$, the only nontrivial group of prime units is $\{1,2\}$ of order $2$. When $r>3$, every group of prime units must have odd order. This is because a subgroup of even order contains $-1$, whose standard representative in $\{2,\dots,r-1\}$ is $r-1$, which is an even integer greater than $2$ and hence composite, which contradicts the definition of a group of prime units. 

A computational search for groups of prime units of orders $3$, $5$, $7$, $9$, and $11$ gives the counts and examples shown in Tables~\ref{tab:prime-group-counts} and~\ref{tab:prime-group-examples}.

\begin{table}[htbp]
\centering
\caption{Number of prime moduli $r\leq X$ having
a group of prime units of order $d$.}
\label{tab:prime-group-counts}
\begin{tabular}{|c|c|c|c|c|c|}
\hline
$X$ & $d=3$ & $d=5$ & $d=7$ & $d=9$ & $d=11$ \\
\hline
$10^{3}$  & 9       & 0   & 0 & 0 & 0 \\
$10^{4}$  & 36      & 0   & 0 & 0 & 0 \\
$10^{5}$  & 188     & 0   & 0 & 0 & 0 \\
$10^{6}$  & 928     & 1   & 0 & 0 & 0 \\
$10^{7}$  & 5462    & 3   & 1 & 0 & 0 \\
$10^{8}$  & 33904   & 8   & 1 & 0 & 0 \\
$10^{9}$  & 227435  & 67  & 1 & 0 & 0 \\
$10^{10}$ & 1606063 & 396 & 4 & 0 & 0 \\
\hline
\end{tabular}
\end{table} 

\begin{table}[htbp]
\centering
\caption{Groups of prime units of orders $3$, $5$, and $7$ modulo primes $r$.}
\label{tab:prime-group-examples}
\small
\renewcommand{\arraystretch}{1.45}
\setlength{\tabcolsep}{8pt}

\begin{tabular}{|lll|}
\hline
\multicolumn{3}{|c|}{$d=3,\quad r\leq10^3$} \\
\hline
$\mathbb{P}_{19,3}^{\times}=\{1,7,11\}$      & $\mathbb{P}_{61,3}^{\times}=\{1,13,47\}$   & $\mathbb{P}_{67,3}^{\times}=\{1,29,37\}$ \\
$\mathbb{P}_{127,3}^{\times}=\{1,19,107\}$   & $\mathbb{P}_{367,3}^{\times}=\{1,83,283\}$ & $\mathbb{P}_{499,3}^{\times}=\{1,139,359\}$ \\
$\mathbb{P}_{571,3}^{\times}=\{1,109,461\}$  & $\mathbb{P}_{631,3}^{\times}=\{1,43,587\}$ & $\mathbb{P}_{991,3}^{\times}=\{1,113,877\}$ \\
\hline
\end{tabular}

\vspace{1.2em}

\begin{tabular}{|l|}
\hline
\multicolumn{1}{|c|}{$d=5,\quad r\leq10^8$} \\
\hline
$\mathbb{P}_{746041,5}^{\times}=\{1,60029,102811,192833,390367\}$ \\
$\mathbb{P}_{1750351,5}^{\times}=\{1,81281,161983,730799,776287\}$ \\
$\mathbb{P}_{7246181,5}^{\times}=\{1,3654143,5534519,5950319,6599561\}$ \\
$\mathbb{P}_{14931271,5}^{\times}=\{1,8505313,11155561,11597617,13535321\}$ \\
$\mathbb{P}_{23445881,5}^{\times}=\{1,9188929,14653549,23095573,23399591\}$ \\
$\mathbb{P}_{24990071,5}^{\times}=\{1,413443,3028609,4693261,16854757\}$ \\
$\mathbb{P}_{28907891,5}^{\times}=\{1,1305401,2550167,9412133,15640189\}$ \\
$\mathbb{P}_{69788531,5}^{\times}=\{1,1099057,3062063,26394701,39232709\}$ \\
\hline
\end{tabular}

\vspace{1.2em}

\begin{tabular}{|l|}
\hline
\multicolumn{1}{|c|}{$d=7,\quad r\leq10^{10}$} \\
\hline
$\mathbb{P}_{2114449,7}^{\times}=\{1,473141,528091,1196773,1294399,1389589,1461353\}$ \\
$\mathbb{P}_{4753987897,7}^{\times}=\{1,240369211,1605309061,1982336591,3323752429,3389967871,3720228527\}$ \\
$\mathbb{P}_{6611641457,7}^{\times}=\{1,341825293,861387953,2128634941,4243083503,5886513439,6373479241\}$ \\
$\mathbb{P}_{9656329789,7}^{\times}=\{1,2040119399,2135930171,2202508507,4624665733,8478112519,9487653037\}$ \\
\hline
\end{tabular}

\end{table}
\FloatBarrier

In the rest of the paper, we focus on the first nontrivial case $d = 3$.

\begin{definition}
A \emph{trey group} is a group of prime units of order $3$ modulo a prime. A \emph{trey prime} is
a prime $r$ for which a trey group exists modulo $r$.\end{definition}

The following theorem gives a criterion for constructing trey groups using three primes.

\begin{theorem}\label{thm:criterion}
Let $p<q<r$ be primes. The following are equivalent:

\begin{enumerate}[label=\textup{(\roman*)}]
    \item $\{1,p,q\}$ is a trey group modulo $r$;
    \item $p^2\equiv q\pmod r$ and $pq\equiv1\pmod r$;
    \item $r=p+q+1$ and $r\mid p^2+p+1$.
\end{enumerate}
\end{theorem}

\begin{proof}
Assume \textup{(i)}. Every nonidentity element of a group of order $3$ has order $3$. Hence $p^2$ is neither $1$ nor $p$. So $p^2\equiv q\pmod r$. The elements $p$ and $q$ are inverses, so $pq\equiv1\pmod r$. This proves \textup{(ii)}.

Assume \textup{(ii)}. We have $p^3\equiv pq\equiv1\pmod r$. Hence $r\mid p^3-1=(p-1)(p^2+p+1)$. Since $p<r$, we have $r\nmid p-1$. Therefore $r\mid p^2+p+1$. Using $p^2\equiv q\pmod r$, we have $p+q+1\equiv0\pmod r$. Since $0<p+q+1<2r$, it follows that $p+q+1=r$. This proves \textup{(iii)}.

Assume \textup{(iii)}. We have $q\equiv-p-1\pmod r$ and $p^2\equiv-p-1\equiv q\pmod r$. It follows that $p^3\equiv p(-p-1)\equiv-p^2-p\equiv1\pmod r$. Thus $p$ has order $3$, and $\{1,p,p^2\}=\{1,p,q\}$ is a subgroup of $\mathbb F_r^\times$. Since $p$ and $q$ are prime, it is a trey group.
\end{proof}

We can also characterize trey primes through the roots of a quadratic congruence. This characterization will be useful when studying how frequently both roots are prime.

\begin{corollary}\label{cor:prime-roots}
A prime $r>3$ is a trey prime if and only if $r\equiv1\pmod3$ and both standard representatives of the roots of $x^2+x+1\equiv0\pmod r$ are prime.
\end{corollary}

\begin{proof}
Since $x^3-1=(x-1)(x^2+x+1)$ and $r>3$, the roots of $x^2+x+1\equiv0\pmod r$ are exactly the elements of order $3$ in $(\mathbb Z/r\mathbb Z)^\times$. Since this group is cyclic of order $r-1$, it has exactly two elements of order $3$ if and only if $3\mid r-1$. These two elements form, together with $1$, the unique subgroup of order $3$ modulo $r$. Therefore, this subgroup is a trey group if and only if the standard representatives of both roots are prime.
\end{proof}

For example, when $r=19$, the roots of $x^2+x+1\equiv0\pmod{19}$ are $7$ and $11$. Since both are prime, $\{1,7,11\}$ is a trey group modulo $19$.

We now prove that there are infinitely many trey primes under Schinzel's Hypothesis H. For this, we use Schinzel's Hypothesis H, which extends Dickson's conjecture from linear polynomials to irreducible polynomials of arbitrary degree.

\begin{conjecture}[Schinzel's Hypothesis H \cite{SchinzelSierpinski}]\label{conj:schinzel}
Let $f_1,\ldots,f_k\in\mathbb Z[x]$ be irreducible polynomials with positive leading coefficients. Suppose no prime divides $f_1(x)\cdots f_k(x)$ for every integer $x$. Then all the polynomials take prime values simultaneously for infinitely many positive integers $x$.
\end{conjecture}

\begin{proof}[Proof of Theorem~\ref{thm:schinzel}]
Assume Schinzel's Hypothesis H. We want to prove that infinitely many different prime moduli have trey groups. By Theorem~\ref{thm:criterion}, it suffices to find infinitely many triples of primes $p<q<r$ such that
$$
r=p+q+1,\qquad r\mid p^2+p+1,
$$
with different values of $r$.

To find such triples, we will construct three polynomials $P(x)$, $Q(x)$, and $R(x)$ whose values satisfy these two relations for every positive integer $x$. We will then use Schinzel's Hypothesis H to prove that all three polynomials take prime values simultaneously for infinitely many positive integers $x$. For every such $x$, if $P(x)<Q(x)<R(x)$, Theorem~\ref{thm:criterion} will give a trey group modulo the prime $R(x)$.

We first construct the three polynomials so that the required relations hold automatically. To ensure that $R(x)$ divides $P(x)^2+P(x)+1$, we can require $P(x)^2+P(x)+1$ to be a fixed positive integer multiple of $R(x)$. We cannot simply set $P(x)^2+P(x)+1=R(x)$, because combining this with the other required relation $R(x)=P(x)+Q(x)+1$ would give $Q(x)=P(x)^2$, which is composite whenever $P(x)>1$. We also cannot set $P(x)^2+P(x)+1=2R(x)$, because $P(x)^2+P(x)+1$ is odd for every integer value of $P(x)$, whereas $2R(x)$ is even. We therefore choose the next possibility and require
$$
P(x)^2+P(x)+1=3R(x).
$$
To ensure the other required relation $R(x)=P(x)+Q(x)+1$, we then define $Q(x)=R(x)-P(x)-1$.

We now need to choose $P(x)$ so that $R(x)=(P(x)^2+P(x)+1)/3$ is a polynomial with integer coefficients. For this, it suffices to choose $P(x)$ congruent to $1$ modulo $3$, since then $P(x)^2+P(x)+1$ is divisible by $3$. A simple choice is $P(x)=3x+1$. However, when $x$ is a positive odd integer, $3x+1$ is an even integer greater than $2$, so it cannot be prime. We therefore restrict to even values of the parameter and choose $P(x)=6x+1$. Substituting this choice into the two relations gives
$$
R(x)=\frac{(6x+1)^2+(6x+1)+1}{3}
=12x^2+6x+1
$$
and
$$
Q(x)=R(x)-P(x)-1=12x^2-1.
$$
Thus we have constructed the three polynomials
$$
P(x)=6x+1,\qquad
Q(x)=12x^2-1,\qquad
R(x)=12x^2+6x+1.
$$
By their construction, they satisfy
$$
P(x)+Q(x)+1=R(x),\qquad
P(x)^2+P(x)+1=3R(x)
$$
for every integer $x$. Therefore, whenever all three polynomial values are prime and satisfy $P(x)<Q(x)<R(x)$, the two arithmetic conditions of Theorem~\ref{thm:criterion} are satisfied.

To prove that $P(x)$, $Q(x)$, and $R(x)$ are simultaneously prime for infinitely many positive integers $x$, we apply Schinzel's Hypothesis H. To use Schinzel's Hypothesis H, we need to verify that all three polynomials are irreducible with positive leading coefficients, and that for every prime $s$, there exists an integer $x$ such that none of the three polynomial values is divisible by $s$.

To verify the first condition, we note that all three polynomials have positive leading coefficients. The polynomial $P(x)=6x+1$ is linear and hence irreducible. To prove that $Q(x)=12x^2-1$ and $R(x)=12x^2+6x+1$ are irreducible, it suffices to show that they have no rational roots. Their discriminants are $48$ and $-12$, respectively. Since neither discriminant is a square in $\mathbb Q$, neither polynomial has a rational root. Therefore, both $Q(x)$ and $R(x)$ are irreducible. Therefore the first condition is satisfied.

To verify the second condition, we must show that no prime $s$ divides $P(x)Q(x)R(x)$ for every integer $x$. It suffices to find one integer $x$ for which the product $P(x)Q(x)R(x)$ is not divisible by any prime $s$. We can take $x=0$, then
$$
P(0)=1,\qquad Q(0)=-1,\qquad R(0)=1,
$$
so $P(0)Q(0)R(0)=-1$. Since no prime divides $-1$, no prime $s$ can divide $P(x)Q(x)R(x)$ for every integer $x$. Thus the second condition is also satisfied.

We have now verified both conditions of Schinzel's Hypothesis H. Therefore, there are infinitely many positive integers $x$ for which all three numbers
$$
p=P(x)=6x+1,\qquad
q=Q(x)=12x^2-1,\qquad
r=R(x)=12x^2+6x+1
$$
are prime.

For every such positive integer $x$, we still need to verify that $p<q<r$ so that Theorem~\ref{thm:criterion} applies. Since $x\geq1$, we have
$$
Q(x)-P(x)=12x^2-6x-2
=6x(2x-1)-2\geq4>0
$$
and
$$
R(x)-Q(x)=6x+2>0.
$$
Therefore, $P(x)<Q(x)<R(x)$. Moreover, since $P(x)+Q(x)+1=R(x)$ and $P(x)^2+P(x)+1=3R(x)$, we have $r=p+q+1$ and $r\mid p^2+p+1$. Thus all the conditions of Theorem~\ref{thm:criterion} are satisfied. Hence $\{1,P(x),Q(x)\}$ is a trey group modulo the prime $R(x)$ for every positive integer $x$ for which all three polynomial values are prime.

Finally, we must show that these trey groups occur modulo infinitely many different prime moduli. Since $R(x)=12x^2+6x+1$ is strictly increasing for positive integers $x$, different positive integers $x$ give different values of $R(x)$. By Schinzel's Hypothesis H, there are infinitely many positive integers $x$ for which $P(x)$, $Q(x)$, and $R(x)$ are simultaneously prime. Therefore, there are infinitely many different prime moduli $R(x)$ having trey groups. This proves the theorem.
\end{proof}

For example, our construction gives the family $T_x=\{1,6x+1,12x^2-1\}$ modulo $r_x=12x^2+6x+1$. By Theorem~\ref{thm:criterion}, this set forms a group of order $3$ modulo $r_x$ whenever the three numbers $6x+1$, $12x^2-1$, and $12x^2+6x+1$ are prime. Taking $x=1$ gives $p=7$, $q=11$, and $r=19$, so $T_1=\{1,7,11\}$ is a trey group modulo $19$. Taking $x=2$ gives $p=13$, $q=47$, and $r=61$, so $T_2=\{1,13,47\}$ is a trey group modulo $61$. Under Schinzel's Hypothesis H, this same polynomial family produces trey groups modulo infinitely many different primes $r_x$.

\section{Sparsity of trey primes}\label{sec:upper}

In this section, we prove the unconditional upper bound $N_3(X)\ll\frac{X}{(\log X)^3}$. We will use the affine sieve. First, we use the relations in a trey group $\{1,p,q\}$ modulo $r$ to construct a binary quadratic form of discriminant $-12$. This allows us to associate the group with a matrix $\gamma\in\mathrm{SL}_2(\mathbb Z)$ satisfying $\|\gamma\|^2<2r$. We then construct a polynomial $F(\gamma)$ such that $F(\gamma)=4rpq$ for these matrices, so that the primality conditions can be studied using a sieve. To apply the sieve, we calculate the local densities of $F$, show that the sieve has dimension $3$, and establish the required congruence distribution of the matrices. These results then give the desired upper bound.

\begin{proof}[Proof of Theorem~\ref{thm:main-upper}]
We want to bound the number of trey primes $r\leq X$. Our approach is to assign to each such prime a matrix $\gamma\in\mathrm{SL}_2(\mathbb Z)$ with $\|\gamma\|<\sqrt{2X}$, and then use the affine sieve to count matrices satisfying the necessary prime conditions, giving an upper bound.

\medskip
\noindent\textit{(i) Matrix representation.}
\par\nobreak\smallskip

\noindent We first construct the matrix representation. Let $\{1,p,q\}$ be a trey group modulo a prime $r$, with $p<q$. To obtain an integer matrix from this group, we will construct a binary quadratic form that is properly equivalent to $X^2+3Y^2$. We begin by looking for a quadratic relation between $r$ and the difference between $p$ and $q$.

Let $a=q-p$. By Theorem~\ref{thm:criterion}, we have
$$
r=p+q+1,\qquad r\mid p^2+p+1.
$$
Since $q=r-p-1$, we have $a=q-p=r-1-2p$, so $a\equiv-(2p+1)\pmod r$. Therefore,
$$
a^2+3\equiv(2p+1)^2+3
=4(p^2+p+1)\equiv0\pmod r.
$$
Thus $r\mid a^2+3$, which gives
$$
rb-a^2=3
$$
for some positive integer $b=(a^2+3)/r$.

We now use this relation to construct the binary quadratic form
$$
f(X,Y)=rX^2+2aXY+bY^2
$$
with discriminant
$$
(2a)^2-4rb=4(a^2-rb)=-12.
$$

The form $f$ is properly equivalent to $X^2+3Y^2$, since the classical classification of binary quadratic forms gives only one proper equivalence class of primitive positive definite forms of discriminant $-12$ \cite{Cox2013}, and $f$ satisfies these conditions: it has discriminant $-12$, it is positive definite because
$$
f(X,Y)=r\left(X+\frac arY\right)^2+\frac3rY^2>0
$$
whenever $(X,Y)\neq(0,0)$, and it is primitive because $r$ is an odd prime and $0<a<r$, so $r\nmid2a$ and hence $(r,2a,b)=1$. In fact, up to proper equivalence, the only primitive reduced positive definite binary quadratic form of discriminant $-12$ is $[1,0,3]$, representing $X^2+3Y^2$.

By this proper equivalence, there exists a matrix
$$
\gamma=
\begin{pmatrix}
x&u\\
y&v
\end{pmatrix}
\in\mathrm{SL}_2(\mathbb Z)
$$
such that
$$
\begin{pmatrix}
r&a\\
a&b
\end{pmatrix}
=
\gamma^{\mathsf T}
\begin{pmatrix}
1&0\\
0&3
\end{pmatrix}
\gamma
=
\begin{pmatrix}
x^2+3y^2&xu+3yv\\
xu+3yv&u^2+3v^2
\end{pmatrix}.
$$

Since $r$ was an arbitrary trey prime, this construction gives a matrix $\gamma\in\mathrm{SL}_2(\mathbb Z)$ for every trey prime $r$. Moreover, if $r\leq X$, this matrix satisfies $\|\gamma\|<\sqrt{2X}$, where $\|\gamma\|=(x^2+y^2+u^2+v^2)^{1/2}$, because $0<a\leq r-2$ implies $b=(a^2+3)/r<r$, and hence
$$
\|\gamma\|^2
=x^2+y^2+u^2+v^2
\leq(x^2+3y^2)+(u^2+3v^2)
=r+b<2r\leq2X.
$$

We have therefore reduced the problem of bounding the number of trey primes to bounding the number of such matrices $\gamma$.

\medskip
\noindent\textit{(ii) The sieve polynomial.}
\par\nobreak\smallskip
However, the number of matrices $\gamma \in \mathrm{SL}_2(\mathbb Z)$ with $\|\gamma\|<\sqrt{2X}$ is of order $X$, which is too big. So the size requirement alone is insufficient to obtain the desired factor $(\log X)^{-3}$. To achieve the target bound, we further use the primality of $p$, $q$, and $r$. The idea is to construct a polynomial $F(\gamma)$ for the matrices $\gamma$ arising from trey groups such that $F(\gamma)=4pqr$ for each such $\gamma$, and use the affine sieve to remove those $\gamma$ whose polynomial value $F(\gamma)$ is divisible by small odd primes. This sieving removes ``bad'' matrices correctly since for ``good'' matrices $\gamma$ arising from trey groups, $F(\gamma)=4pqr$ is only divisible by $2$, $p$, $q$, and $r$, and it should not be divisible by small odd primes such as those smaller than $p$.

To construct this polynomial, we extend the expressions for $r$ and $a$ to any matrix
$$
\gamma=
\begin{pmatrix}
x&u\\
y&v
\end{pmatrix}
\in\mathrm{SL}_2(\mathbb Z)
$$
by defining
$$
R(\gamma)=x^2+3y^2,\qquad
A(\gamma)=xu+3yv.
$$
For a matrix arising from a trey group $\{1,p,q\}$ modulo $r$, we have $R(\gamma)=r$ and $A(\gamma)=q-p$. Since $r = p+q+1$, we have
$$
2p=R(\gamma)-1-A(\gamma),\qquad
2q=R(\gamma)-1+A(\gamma).
$$
We can then define the polynomial as
$$
F(\gamma)=R(\gamma)
\bigl(R(\gamma)-1-A(\gamma)\bigr)
\bigl(R(\gamma)-1+A(\gamma)\bigr),
$$
which satisfies
\begin{equation}\label{eq:F-trey}
F(\gamma)=4pqr
\end{equation}
for every matrix $\gamma$ arising from a trey group.

This identity allows us to use the sieve. In particular, when $p$, $q$, and $r$ are sufficiently large, the integer
$F(\gamma) = 4pqr$ has no small odd prime divisors.

For a prime $\ell>3$, define the local density at $\ell$ to be
$$
\beta(\ell)=
\frac{\#\{\gamma\in\mathrm{SL}_2(\mathbb F_\ell):
F(\gamma)=0\}}
{\#\mathrm{SL}_2(\mathbb F_\ell)}.
$$
We will show that
$$
\beta(\ell)
=\frac{\ell(3+\chi_\ell)-3}{\ell(\ell+\chi_\ell)},
$$
where $\chi_\ell=\left(\frac{-3}{\ell}\right)$ is the Legendre symbol.

\medskip
\noindent\textit{(iii) Local densities.}
\par\nobreak\smallskip
To calculate $\beta(\ell)$, we need to count matrices $\gamma\in\mathrm{SL}_2(\mathbb F_\ell)$ for which $F(\gamma)=0$. Since
$$
F(\gamma)=R(\gamma)
\bigl(R(\gamma)-1-A(\gamma)\bigr)
\bigl(R(\gamma)-1+A(\gamma)\bigr),
$$
this means that $R=0$, $R-1-A=0$, or $R-1+A=0$. Instead of counting the four matrix entries directly, we can count these conditions using a simpler equation involving $R$ and $A$. For this, we define
$$
B(\gamma)=u^2+3v^2,
$$
we will have
\begin{equation}\label{eq:determinant-identity}
R(\gamma)B(\gamma)-A(\gamma)^2
=3(xv-yu)^2=3,
\end{equation}
since $\gamma\in\mathrm{SL}_2(\mathbb F_\ell)$. This way we have extended the relation $rb-a^2=3$, originally obtained from a trey group, to the identity $R(\gamma)B(\gamma)-A(\gamma)^2=3$ for every matrix $\gamma\in\mathrm{SL}_2(\mathbb F_\ell)$.

Thus every matrix gives a triple $(R,A,B)$ satisfying $RB-A^2=3$. If every such triple comes from the same number of matrices, we can calculate $\beta(\ell)$ by counting triples instead of matrices, since both the numerator and denominator of $\beta(\ell)$ are multiplied by this same number. We now prove that they indeed come from the same number of matrices.

More precisely, let
$$
V_\ell=\{(R,A,B)\in\mathbb F_\ell^3:RB-A^2=3\},
$$
and define
$$
\Phi_\ell:\mathrm{SL}_2(\mathbb F_\ell)\longrightarrow V_\ell,
\qquad
\Phi_\ell(\gamma)=(R(\gamma),A(\gamma),B(\gamma)).
$$
We need to show that $\Phi_\ell$ is onto and that every triple in $V_\ell$ has the same number of preimages.

First, we prove that $\Phi_\ell$ is onto. Fix any $(R,A,B)\in V_\ell$. We need to find $x,y,u,v\in\mathbb F_\ell$ satisfying
$$
x^2+3y^2=R,\qquad
xu+3yv=A,\qquad
u^2+3v^2=B,\qquad
xv-yu=1.
$$

Suppose first that $R\neq0$. We begin by finding $x,y$ such that $x^2+3y^2=R$. The two sets
$$
\{x^2:x\in\mathbb F_\ell\}
\quad\text{and}\quad
\{R-3y^2:y\in\mathbb F_\ell\}
$$
each contain $(\ell+1)/2$ elements. Since their sizes add up to $\ell+1$, they must intersect in $\mathbb F_\ell$. Hence we can choose $(x,y)$ satisfying $x^2+3y^2=R$.

With this $(x,y)$, the equations
$$
xu+3yv=A,\qquad xv-yu=1
$$
form a linear system in $u,v$ whose determinant is $x^2+3y^2=R\neq0$. Thus they have the solution
$$
u=\frac{Ax-3y}{R},\qquad
v=\frac{Ay+x}{R}.
$$
The remaining equation also holds, since
$$
u^2+3v^2=\frac{A^2+3}{R}=B,
$$
where the last equality follows from $RB-A^2=3$.

Now suppose $R=0$. Since $RB-A^2=3$, we have $A^2=-3$ in $\mathbb F_\ell$, so $A\neq0$. To satisfy $x^2+3y^2=0$, we choose $x=-A$ and $y=1$. The equation $xv-yu=1$ then gives $u=-Av-1$. Substituting this into $u^2+3v^2=B$ gives $2Av+1=B$, so we choose
$$
v=\frac{B-1}{2A},\qquad u=-Av-1.
$$
These values also satisfy
$$
xu+3yv=A+(A^2+3)v=A.
$$
Thus all four required equations hold.

Therefore, every triple in $V_\ell$ comes from
at least one matrix $\gamma=\left(\begin{smallmatrix}x&u\\y&v\end{smallmatrix}\right) \in \mathrm{SL}_2(\mathbb F_\ell)$, proving
that $\Phi_\ell$ is onto.

Next, we prove that every triple in $V_\ell$ comes from the same number of matrices. Fix a matrix $\gamma_0$ giving a triple $(R,A,B)$. We want to find all matrices $\gamma$ giving this same triple.

To compare these matrices, we write $R$, $A$, and $B$ in a single matrix equation. Let
$$
D=\begin{pmatrix}1&0\\0&3\end{pmatrix}.
$$
We have
$$
\begin{aligned}
\gamma^{\mathsf T}D\gamma
&=
\begin{pmatrix}x&y\\u&v\end{pmatrix}
\begin{pmatrix}1&0\\0&3\end{pmatrix}
\begin{pmatrix}x&u\\y&v\end{pmatrix} =
\begin{pmatrix}
x^2+3y^2&xu+3yv\\
xu+3yv&u^2+3v^2
\end{pmatrix} =
\begin{pmatrix}
R(\gamma)&A(\gamma)\\
A(\gamma)&B(\gamma)
\end{pmatrix}.
\end{aligned}
$$
The same calculation holds for $\gamma_0$. Thus $\gamma$ and $\gamma_0$ give the same triple if and only if
$$
\gamma^{\mathsf T}D\gamma
=\gamma_0^{\mathsf T}D\gamma_0.
$$

To find all matrices satisfying this equation, write $\gamma=h\gamma_0$, where $h=\gamma\gamma_0^{-1}\in\mathrm{SL}_2(\mathbb F_\ell)$. Then
$$
\begin{aligned}
\gamma^{\mathsf T}D\gamma
=\gamma_0^{\mathsf T}D\gamma_0
&\iff
\gamma_0^{\mathsf T}h^{\mathsf T}Dh\gamma_0
=\gamma_0^{\mathsf T}D\gamma_0 \iff h^{\mathsf T}Dh=D.
\end{aligned}
$$

We therefore define
$$
H_\ell=
\{h\in\mathrm{SL}_2(\mathbb F_\ell):h^{\mathsf T}Dh=D\}.
$$
The matrices giving the same triple as $\gamma_0$ are exactly those of the form $\gamma=h\gamma_0$ with $h\in H_\ell$. Since $\gamma_0$ is invertible, $h_1\gamma_0=h_2\gamma_0$ implies $h_1=h_2$. Thus different elements of $H_\ell$ give different matrices $\gamma$, so exactly $\#H_\ell$ matrices give the same triple as $\gamma_0$. This number does not depend on the chosen triple, so every triple in $V_\ell$ comes from the same number of matrices.

We have now proved that $\Phi_\ell$ is onto and that all its fibers have the same size. Therefore, the proportion of matrices for which $F(\gamma)=0$ equals the proportion of triples in $V_\ell$ satisfying $R(R-1-A)(R-1+A)=0$. That is,
$$
\beta(\ell)=
\frac{
\#\{(R,A,B)\in V_\ell:
R(R-1-A)(R-1+A)=0\}
}{
\#V_\ell
}.
$$

To get $\beta(\ell)$, we first calculate $\#V_\ell$. When $R\neq0$, the equation $RB-A^2=3$ uniquely determines
$$
B=\frac{A^2+3}{R}
$$
for every choice of $R\neq0$ and $A$. This gives $\ell(\ell-1)$ points.

When $R=0$, the equation becomes $A^2=-3$. By the definition of the Legendre symbol, this equation has two solutions if $\chi_\ell=1$ and no solutions if $\chi_\ell=-1$. Thus there are $1+\chi_\ell$ possible values of $A$. Since $B$ is unrestricted, each value of $A$ allows $\ell$ choices of $B$. Therefore,
$$
\#V_\ell
=\ell(\ell-1)+\ell(1+\chi_\ell)
=\ell(\ell+\chi_\ell).
$$

Next, we count the points where at least one of the three factors $R$, $R-1-A$, and $R-1+A$ vanishes.

The condition $R=0$ gives $\ell(1+\chi_\ell)$ points, as calculated above.

The condition $R-1-A=0$ gives $A=R-1$. If $R=0$, then $A=-1$, and the equation $RB-A^2=3$ would give $-1=3$ in $\mathbb F_\ell$, which is impossible since $\ell>3$. Thus $R\neq0$. For each of the $\ell-1$ nonzero values of $R$, the value $A=R-1$ is fixed, and
$$
B=\frac{A^2+3}{R}
$$
is uniquely determined. Hence this condition gives $\ell-1$ points.

The same argument shows that $R-1+A=0$ also gives $\ell-1$ points.

The first set is disjoint from the other two. In fact, $R=0$ together with either $R-1-A=0$ or $R-1+A=0$ would give $A=\pm1$, while $RB-A^2=3$ would require $A^2=-3$. This would imply $4=0$ in $\mathbb F_\ell$, which is impossible for $\ell>3$.

The last two sets intersect only when $R=1$ and $A=0$. The equation $RB-A^2=3$ then gives $B=3$, so their intersection consists of the single point $(1,0,3)$.

Subtracting this overlap, we have
$$
\begin{aligned}
&\#\{(R,A,B)\in V_\ell:
R(R-1-A)(R-1+A)=0\}\\
&\qquad
=\ell(1+\chi_\ell)+2(\ell-1)-1\\
&\qquad
=\ell(3+\chi_\ell)-3.
\end{aligned}
$$

Since all fibers of $\Phi_\ell$ have equal size, the local density is therefore
$$
\beta(\ell)
=\frac{\ell(3+\chi_\ell)-3}{\ell(\ell+\chi_\ell)}
=\frac{3+\chi_\ell}{\ell}+O(\ell^{-2}).
$$
In particular, $\beta(\ell)>0$.

To apply the sieve and obtain a logarithmic saving of order $(\log z)^{-3}$ for some small number $z$, we need to establish two properties of the local densities $\beta(\ell)$. First, $1-\beta(\ell)$ must have a positive lower bound independent of $\ell$. Second, we need to verify that the sieve has dimension $3$ for the factor $(\log X)^{-3}$ in the upper bound.

We first show that $\beta(\ell)\leq1/2$ for every prime $\ell>3$. Recall that
$$
\beta(\ell)=
\frac{\ell(3+\chi_\ell)-3}{\ell(\ell+\chi_\ell)},
\qquad
\chi_\ell=\left(\frac{-3}{\ell}\right).
$$
If $\chi_\ell=-1$, then
$$
\beta(\ell)=\frac{2\ell-3}{\ell(\ell-1)},
$$
and the inequality $\beta(\ell)\leq1/2$ is equivalent to
$$
(\ell-2)(\ell-3)\geq0,
$$
which holds since $\ell\geq5$.

If $\chi_\ell=1$, then
$$
\beta(\ell)=\frac{4\ell-3}{\ell(\ell+1)},
$$
and the inequality $\beta(\ell)\leq1/2$ is equivalent to
$$
(\ell-1)(\ell-6)\geq0,
$$
which holds since $\ell\geq7$. Therefore,
$$
0<\beta(\ell)\leq\frac12
$$
for every prime $\ell>3$, so $1-\beta(\ell)\geq1/2$ for every prime $\ell > 3$.

\medskip
\noindent\textit{(iv) Sieve dimension.}
\par\nobreak\smallskip
We now verify that the sieve has dimension $3$ to obtain the factor $(\log X)^{-3}$ in the upper bound. Following the sieve-dimension condition $(A_2)$ of Bourgain, Gamburd, and Sarnak \cite{BGS10}, for the local densities $\beta(\ell)$, the sieve dimension $\kappa$ is characterized by
\begin{equation}\label{eq:sieve-dimension-condition}
\sum_{\substack{w\leq\ell<z\\
\ell\ {\rm prime}}}
\beta(\ell)\log\ell
=\kappa\log(z/w)+O(1)
\qquad(3<w<z).
\end{equation}

This $\kappa$ determines the exponent of $\log z$ in the estimate effect
$$
\prod_{3<\ell<z}(1-\beta(\ell))
\asymp(\log z)^{-\kappa},
$$ which is the combined effect of
excluding matrices for which $F(\gamma)\equiv0\pmod\ell$ for primes $3<\ell<z$. Proving $\kappa = 3$ is directly associated to the factor $(\log X)^{-3}$ in our final goal.

We therefore need to prove that the coefficient $\kappa$ of $\log(z/w)$ in the sum above is $3$.

Recall that
$$
\beta(\ell)=
\frac{\ell(3+\chi_\ell)-3}
{\ell(\ell+\chi_\ell)},
\qquad
\chi_\ell=\left(\frac{-3}{\ell}\right).
$$
To determine the coefficient of $\log(z/w)$, we first separate the leading term of $\beta(\ell)$:
$$
\begin{aligned}
\beta(\ell)
=\frac{\ell(3+\chi_\ell)-3}
{\ell(\ell+\chi_\ell)} 
=\frac{3+\chi_\ell}{\ell}
-\frac{4+3\chi_\ell}{\ell(\ell+\chi_\ell)}.
\end{aligned}
$$

By quadratic reciprocity,
$$
\chi_\ell
=\left(\frac{-3}{\ell}\right)
=\left(\frac{\ell}{3}\right)
=
\begin{cases}
1,&\ell\equiv1\pmod3,\\
-1,&\ell\equiv2\pmod3.
\end{cases}
$$
Hence
$$
\beta(\ell)=
\begin{cases}
\displaystyle
\frac4\ell-\frac7{\ell(\ell+1)},
&\ell\equiv1\pmod3,\\[8pt]
\displaystyle
\frac2\ell-\frac1{\ell(\ell-1)},
&\ell\equiv2\pmod3.
\end{cases}
$$
Therefore,
$$
\beta(\ell)=
\begin{cases}
4/\ell+O(\ell^{-2}),
&\ell\equiv1\pmod3,\\
2/\ell+O(\ell^{-2}),
&\ell\equiv2\pmod3.
\end{cases}
$$

Plug this in the sum we have
\begin{equation}\label{eq:beta-sum-expansion}
\begin{aligned}
\sum_{3<\ell<z}\beta(\ell)\log\ell
&=
4\sum_{\substack{3<\ell<z\\
\ell\equiv1\pmod3}}\frac{\log\ell}{\ell} +
2\sum_{\substack{3<\ell<z\\
\ell\equiv2\pmod3}}\frac{\log\ell}{\ell} +
O\left(\sum_{\ell>3}\frac{\log\ell}{\ell^2}\right).
\end{aligned}
\end{equation}
Note that the error term is bounded independently of $z$, because
$$
\begin{aligned}
\sum_{\substack{\ell>3\\ \ell\ {\rm prime}}}
\frac{\log\ell}{\ell^2}
&\leq
\sum_{n=5}^{\infty}\frac{\log n}{n^2} \leq
\int_4^\infty\frac{\log t}{t^2}\,dt =\frac{\log4+1}{4}<\infty.
\end{aligned}
$$ Hence what remains to do here is to calculate
the two sums of $(\log\ell)/\ell$ over primes congruent to $1$ and $2$ modulo $3$.

To calculate the two sums of $(\log\ell)/\ell$, we use the prime number theorem in arithmetic progressions. The theorem provides estimates for sums of $\log\ell$ over primes in each residue class modulo $3$. We therefore define for $j = 1,2$ that
$$
\Theta_j(t)=
\sum_{\substack{3<\ell\leq t\\
\ell\equiv j\pmod3}}\log\ell.
$$
Since $j = 1,2$ are the two reduced residue classes modulo $3$, the prime number theorem in arithmetic progressions for the fixed modulus $3$ \cite{Davenport2000} gives
$$
\Theta_j(t)
=\frac t2+
O\left(t e^{-c\sqrt{\log t}}\right),
\qquad j\in\{1,2\},
$$
for some constant $c>0$.

We now need to introduce the factor $1/\ell$ into the sums. For every $\ell\leq t$, we have
$$
\begin{aligned}
\frac1\ell
&=\frac1t+\left[-\frac1\tau\right]_\ell^t=\frac1t+\int_\ell^t\frac{d\tau}{\tau^2}.
\end{aligned}
$$
Multiplying by $\log\ell$ and summing over primes in either residue class, we have, for $j=1,2$, that
$$
\begin{aligned}
\sum_{\substack{3<\ell\leq t\\
\ell\equiv j\pmod3}}
\frac{\log\ell}{\ell}
&=
\sum_{\substack{3<\ell\leq t\\
\ell\equiv j\pmod3}}
\log\ell
\left(
\frac1t+\int_\ell^t\frac{d\tau}{\tau^2}
\right)\\
&=
\frac1t
\sum_{\substack{3<\ell\leq t\\
\ell\equiv j\pmod3}}\log\ell \quad +
\sum_{\substack{3<\ell\leq t\\
\ell\equiv j\pmod3}}
\log\ell
\left(
\int_\ell^t\frac{d\tau}{\tau^2}
\right)\\
&=
\frac{\Theta_j(t)}t+
\int_3^t
\left(
\sum_{\substack{3<\ell\leq\tau\\
\ell\equiv j\pmod3}}\log\ell
\right)\frac{d\tau}{\tau^2}\\
&=
\frac{\Theta_j(t)}t+
\int_3^t\frac{\Theta_j(\tau)}{\tau^2}\,d\tau,
\end{aligned}
$$
where the third equality follows by interchanging the finite sum and the integral. This interchange is valid because the sum is finite and that for each fixed $\tau$ with $3\leq\tau\leq t$, the original integration range $\ell\leq\tau\leq t$ restricts the inner sum to primes satisfying $3<\ell\leq\tau$ and $\ell\equiv j\pmod3$.

We next substitute the estimate for $\Theta_j$ into the two terms. For the first term, we have
$$
\frac{\Theta_j(t)}t
=
\frac{t/2+O\left(t e^{-c\sqrt{\log t}}\right)}t
=
\frac12+O\left(e^{-c\sqrt{\log t}}\right).
$$
For the integral term, we have
$$
\begin{aligned}
\int_3^t\frac{\Theta_j(\tau)}{\tau^2}\,d\tau
&=
\int_3^t
\frac{\tau/2+
O\left(\tau e^{-c\sqrt{\log\tau}}\right)}
{\tau^2}\,d\tau\\
&=
\frac12\int_3^t\frac{d\tau}{\tau} +
O\left(
\int_3^t
\frac{e^{-c\sqrt{\log\tau}}}{\tau}\,d\tau
\right)\\
&=
\frac12\log(t/3) +
O\left(
\int_3^t
\frac{e^{-c\sqrt{\log\tau}}}{\tau}\,d\tau
\right).
\end{aligned}
$$

Combining the estimates for the two terms, we have
$$
\begin{aligned}
\sum_{\substack{3<\ell\leq t\\
\ell\equiv j\pmod3}}
\frac{\log\ell}{\ell}
&=
\frac{\Theta_j(t)}t
+\int_3^t\frac{\Theta_j(\tau)}{\tau^2}\,d\tau\\
&=
\frac12+\frac12\log(t/3)\\
&\quad+
O\left(
e^{-c\sqrt{\log t}}
+\int_3^t
\frac{e^{-c\sqrt{\log\tau}}}{\tau}\,d\tau
\right)\\
&=
\frac12\log t+\frac{1-\log3}{2}
+
O\left(
e^{-c\sqrt{\log t}}
+\int_3^t
\frac{e^{-c\sqrt{\log\tau}}}{\tau}\,d\tau
\right)\\
&=
\frac12\log t+
O\left(
1+\int_3^t
\frac{e^{-c\sqrt{\log\tau}}}{\tau}\,d\tau
\right).
\end{aligned}
$$
Here $(1-\log3)/2$ is a constant, and $e^{-c\sqrt{\log t}}\leq1$ for $t\geq3$. Both terms are therefore absorbed into $O(1)$.

We now return to \eqref{eq:beta-sum-expansion} and use the formula just obtained for $j=1$ and $j=2$. We have
$$
\begin{aligned}
\sum_{3<\ell<z}\beta(\ell)\log\ell
&=
4\left(
\frac12\log z+
O\left(
1+\int_3^z
\frac{e^{-c\sqrt{\log\tau}}}{\tau}\,d\tau
\right)
\right)\\
&\quad+
2\left(
\frac12\log z+
O\left(
1+\int_3^z
\frac{e^{-c\sqrt{\log\tau}}}{\tau}\,d\tau
\right)
\right)
+O(1)\\
&=
3\log z+
O\left(
1+\int_3^z
\frac{e^{-c\sqrt{\log\tau}}}{\tau}\,d\tau
\right).
\end{aligned}
$$
Here we have taken $t=z$. Note that there is a slight difference between the formula just obtained with $t=z$ and \eqref{eq:beta-sum-expansion}: the former uses $\ell\leq z$, while the latter uses $\ell<z$, but the difference is at most $(\log z)/z=O(1)$, so it can be absorbed into the error term.

We now return to \eqref{eq:sieve-dimension-condition}. To verify this relation with $\kappa=3$, it is enough to prove that, for every $z>3$, $\sum_{3<\ell<z}\beta(\ell)\log\ell =3\log z+O(1)$ because
$$
\begin{aligned}
\sum_{\substack{w\leq\ell<z\\
\ell\ {\rm prime}}}
\beta(\ell)\log\ell
&=
\sum_{3<\ell<z}\beta(\ell)\log\ell
-
\sum_{3<\ell<w}\beta(\ell)\log\ell\\
&=
\bigl(3\log z+O(1)\bigr)
-
\bigl(3\log w+O(1)\bigr)\\
&=3\log(z/w)+O(1),
\end{aligned}
$$
as required by \eqref{eq:sieve-dimension-condition}.

The previous calculation has already given the main term $3\log z$. To establish the estimate with an $O(1)$ error, it remains to prove that $\int_3^z \frac{e^{-c\sqrt{\log\tau}}}{\tau}\,d\tau =O(1).$

Since the integrand is positive, we have
$$
0\leq
\int_3^z
\frac{e^{-c\sqrt{\log\tau}}}{\tau}\,d\tau
\leq
\int_3^\infty
\frac{e^{-c\sqrt{\log\tau}}}{\tau}\,d\tau.
$$
To evaluate the integral on the right, let $\rho=\sqrt{\log\tau}$. Then $\tau=e^{\rho^2}$ and $d\tau/\tau=2\rho\,d\rho$. Hence
$$
\begin{aligned}
\int_3^\infty
\frac{e^{-c\sqrt{\log\tau}}}{\tau}\,d\tau
&=
2\int_{\sqrt{\log3}}^\infty
\rho e^{-c\rho}\,d\rho =
2\left[
-e^{-c\rho}
\left(\frac{\rho}{c}+\frac1{c^2}\right)
\right]_{\sqrt{\log3}}^\infty =
2e^{-c\sqrt{\log3}}
\left(
\frac{\sqrt{\log3}}c+\frac1{c^2}
\right) <\infty.
\end{aligned}
$$

Since this finite upper bound depends only on the fixed constant $c$ and not on $z$, we have $\int_3^z \frac{e^{-c\sqrt{\log\tau}}}{\tau}\,d\tau =O(1).$ Therefore, the error integral is $O(1)$, and hence
$$
\sum_{3<\ell<z}\beta(\ell)\log\ell
=3\log z+O(1).
$$
Applying this estimate at $z$ and $w$ and subtracting, as mentioned earlier, we have \eqref{eq:sieve-dimension-condition} with $\kappa=3$. Thus the sieve has dimension $3$.

Note that in the sieve application below, the prime divisors of a fixed integer $q_0$ are excluded. Their omission changes the weighted prime sums by $O_{q_0}(1)$, so the sieve-dimension condition remains valid with $\kappa=3$.

\medskip
\noindent\textit{(v) The local-factor product.}
\par\nobreak\smallskip
After proving that the sieve dimension is $3$, we now turn to the local-factor product introduced earlier. Remind that our final goal is to prove $N_3(X)\ll\frac{X}{(\log X)^3}.$ To obtain the logarithmic factor in this upper bound, we need to prove
$$
\prod_{3<\ell<z}(1-\beta(\ell))
\asymp(\log z)^{-3},
$$
where we recall that the product
$$
\prod_{3<\ell<z}(1-\beta(\ell))
=
\prod_{3<\ell<z}
\frac{
\#\{\gamma\in\mathrm{SL}_2(\mathbb F_\ell):
F(\gamma)\not\equiv0\pmod\ell\}
}{
\#\mathrm{SL}_2(\mathbb F_\ell)
}
$$
combines the local proportions of matrices $\gamma$ satisfying $F(\gamma)\not\equiv0\pmod\ell$ for every prime $3<\ell<z$. We will derive its estimate $(\log z)^{-3}$ from the sum of $\beta(\ell)\log\ell$ we have just calculated above.

To estimate the product, we first take its logarithm. Since $0<\beta(\ell)\leq1/2 <1$, we use the Taylor expansion of $\log(1-x)$ at $x=0$:
$$
\log(1-x)
=-x-\sum_{n=2}^{\infty}\frac{x^n}{n},
\qquad 0<x<1,
$$ substituting $x=\beta(\ell)$, we have
$$
\begin{aligned}
\log(1-\beta(\ell))
&=-\beta(\ell)
-\sum_{n=2}^{\infty}\frac{\beta(\ell)^n}{n} =-\beta(\ell)+O(\beta(\ell)^2),
\end{aligned}
$$ where the last equality follows from
$0<\beta(\ell)\leq1/2$, which implies
$$
\begin{aligned}
0
&\leq
\sum_{n=2}^{\infty}\frac{\beta(\ell)^n}{n}
&&\text{(since $\beta(\ell)>0$)}\\
&\leq
\frac12\sum_{n=2}^{\infty}\beta(\ell)^n
&&\text{(since $1/n\leq1/2$)}\\
&=
\frac{\beta(\ell)^2}{2}
\sum_{n=0}^{\infty}\beta(\ell)^n
&&\text{(factor out $\beta(\ell)^2$ and reindex)}\\
&=
\frac{\beta(\ell)^2}{2(1-\beta(\ell))}
&&\text{(geometric series formula)}\\
&\leq
\beta(\ell)^2
&&\text{(since $\beta(\ell)\leq1/2$)}.
\end{aligned}
$$

Summing the Taylor expansion over primes $3<\ell<z$, we have
$$
\begin{aligned}
\log\prod_{3<\ell<z}(1-\beta(\ell))
=
\sum_{3<\ell<z}\log(1-\beta(\ell)) =
-\sum_{3<\ell<z}\beta(\ell)
+
O\left(
\sum_{3<\ell<z}\beta(\ell)^2
\right).
\end{aligned}
$$

To show that $\prod_{3<\ell<z}(1-\beta(\ell)) \asymp(\log z)^{-3},$ we need to show that
$$
\log\prod_{3<\ell<z}(1-\beta(\ell))
=-3\log\log z+O(1).
$$ In fact, the error term
$O\left(\sum_{3<\ell<z}\beta(\ell)^2\right)=O(1).$ This is because
$$
\begin{aligned}
0&\leq
\sum_{3<\ell<z}\beta(\ell)^2 &&\text{(since $\beta(\ell)>0$)}\\
&\leq
\sum_{\substack{\ell>3 \\
\ell\ {\rm prime}}}\beta(\ell)^2
&&\text{(nonnegative terms)}\\
&\leq
16\sum_{\substack{\ell>3\\
\ell\ {\rm prime}}}\frac1{\ell^2}
&&\begin{gathered}
\text{(since $\beta(\ell)=4/\ell-7/(\ell(\ell+1))$}\\
\text{or $2/\ell-1/(\ell(\ell-1))$,}\\
\text{so $0<\beta(\ell)\leq4/\ell$)}
\end{gathered}\\
&\leq
16\sum_{n=5}^{\infty}\frac1{n^2}
&&\text{(sum over all integers $n\geq5$)}\\
&<\infty
&&\text{(convergence of the $p$-series, $p=2$)}.
\end{aligned}
$$
Note that the finite upper bound is independent of $z$. Hence the error term is $O(1)$. Therefore,
$$
\log\prod_{3<\ell<z}(1-\beta(\ell))
=
-\sum_{3<\ell<z}\beta(\ell)+O(1).
$$
To determine the product, it remains to estimate $\sum_{3<\ell<z}\beta(\ell).$ We want to show that
$$
\sum_{3<\ell<z}\beta(\ell)
=3\log\log z+O(1).
$$

Recall that when proving that the sieve dimension is $3$, we obtained $\sum_{3<\ell<t}\beta(\ell)\log\ell =3\log t+O(1).$ Now let $S(t)=\sum_{3<\ell<t}\beta(\ell)\log\ell.$ We have $S(t)=3\log t+O(1).$ For $3<\ell<z$, since $\frac{d}{dt}\left(\frac1{\log t}\right) =-\frac1{t(\log t)^2},$ we have
$$
\begin{aligned}
\int_\ell^z\frac{dt}{t(\log t)^2}
&=\left[-\frac1{\log t}\right]_\ell^z =\frac1{\log\ell}-\frac1{\log z}.
\end{aligned}
$$
Hence $\frac1{\log\ell} =\frac1{\log z} +\int_\ell^z\frac{dt}{t(\log t)^2}.$ Multiplying by $\beta(\ell)\log\ell$ and summing over $3<\ell<z$, we have
$$
\begin{aligned}
\sum_{3<\ell<z}\beta(\ell)
=
\frac1{\log z}
\sum_{3<\ell<z}\beta(\ell)\log\ell +
\sum_{3<\ell<z}\beta(\ell)\log\ell
\int_\ell^z\frac{dt}{t(\log t)^2}=
\frac{S(z)}{\log z}
+\int_3^z\frac{S(t)}{t(\log t)^2}\,dt,
\end{aligned}
$$
where the last equality follows by interchanging the finite sum and the integral.

Substituting $S(t)=3\log t+O(1)$ gives
$$
\begin{aligned}
\sum_{3<\ell<z}\beta(\ell)
=
\frac{3\log z+O(1)}{\log z}
+\int_3^z
\frac{3\log t+O(1)}
{t(\log t)^2}\,dt =
3+3\int_3^z\frac{dt}{t\log t} +
O\left(
1+\int_3^z
\frac{dt}{t(\log t)^2}
\right).
\end{aligned}
$$

The main integral is
$$
\begin{aligned}
\int_3^z\frac{dt}{t\log t}
&=
\left[\log\log t\right]_3^z =
\log\log z-\log\log3.
\end{aligned}
$$
The error integral is
$$
\begin{aligned}
0\leq
\int_3^z\frac{dt}{t(\log t)^2}
=
\left[-\frac1{\log t}\right]_3^z =
\frac1{\log3}-\frac1{\log z} \leq
\frac1{\log3}.
\end{aligned}
$$
Thus the error integral is $O(1)$, since $1/\log3$ is independent of $z$.

Substituting these two integrals back gives
$$
\begin{aligned}
\sum_{3<\ell<z}\beta(\ell)
=
3+3(\log\log z-\log\log3)+O(1) =
3\log\log z+O(1),
\end{aligned}
$$
where the constant $3-3\log\log3$ is absorbed into $O(1)$.

Substituting this estimate into the logarithmic product formula, we have
$$
\begin{aligned}
\log\prod_{3<\ell<z}(1-\beta(\ell))
=
-\sum_{3<\ell<z}\beta(\ell)+O(1) =
-3\log\log z+O(1).
\end{aligned}
$$
Exponentiating both sides gives
$$
\begin{aligned}
\prod_{3<\ell<z}(1-\beta(\ell))
=
\exp\bigl(-3\log\log z+O(1)\bigr) =
(\log z)^{-3}\exp(O(1)) \asymp(\log z)^{-3}.
\end{aligned}
$$
The last relation holds because $\exp(O(1))$ is bounded above and below by positive constants independent of $z$. Thus the local-factor product estimate stated earlier is proved.

More generally, if $S$ is any fixed finite set of primes greater than $3$, then $\prod_{3<\ell<z,\ \ell\notin S}(1-\beta(\ell))\asymp_S(\log z)^{-3}$, because the omitted factors form a finite product of nonzero constants. In the sieve below we take $S$ to be the set of prime divisors of $q_0$ greater than $3$.

\medskip
\noindent\textit{(vi) Congruence distribution.}
\par\nobreak\smallskip
We have just calculated the local densities, but to apply the sieve we also need to know that the integer matrices $\gamma\in\mathrm{SL}_2(\mathbb Z)$ are sufficiently well distributed among residue classes modulo the sieve moduli, otherwise the local densities $\beta(\ell)$ computed over finite fields cannot be used to estimate the corresponding divisibility frequencies among integer matrices, and the sieve error terms cannot be controlled.

We now establish the congruence distribution needed for the sieve. We will show that the proportion of matrices $\gamma$ for which $F(\gamma)$ is divisible by an integer $m$ is approximately the product of the local densities with a controlled error.

Write
$$
M(T)=\#\{\gamma\in\mathrm{SL}_2(\mathbb Z):\|\gamma\|\leq T\}.
$$
We use the congruence-counting theorem of Bourgain--Gamburd--Sarnak \cite[Theorem 1.3]{BGS11}. This theorem applies to $\Lambda=\mathrm{SL}_2(\mathbb Z)$ because the group is finitely generated, has critical exponent $\delta(\Lambda) = 1>1/2$, and reduces surjectively onto $\mathrm{SL}_2(\mathbb Z/m\mathbb Z)$. The cited estimate applies uniformly to squarefree $m$ coprime to a fixed integer $q_0$.

By the theorem, there exist constants $c_0>0$, $\varepsilon>0$, and an integer $q_0\geq1$ such that, for every squarefree $m$ with $(m,q_0)=1$ and every $g\in\mathrm{SL}_2(\mathbb Z/m\mathbb Z)$,
\begin{equation}\label{eq:congruence-count}
\begin{aligned}
&\#\{\gamma\in\mathrm{SL}_2(\mathbb Z):
\|\gamma\|\leq T,\ \gamma\equiv g\pmod m\}\\
&\qquad=
\frac{c_0T^2}{\#\mathrm{SL}_2(\mathbb Z/m\mathbb Z)}
+O(m^3T^{2-\varepsilon}),
\end{aligned}
\end{equation}
where the implied constant in the error term is independent of $m$, $g$, and $T$.

Taking $m=1$ gives
$$
M(T)=c_0T^2+O(T^{2-\varepsilon})\asymp T^2.
$$
Thus \eqref{eq:congruence-count} can also be written as
$$
\begin{aligned}
&\#\{\gamma\in\mathrm{SL}_2(\mathbb Z):
\|\gamma\|\leq T,\ \gamma\equiv g\pmod m\}\\
&\qquad=
\frac{M(T)}{\#\mathrm{SL}_2(\mathbb Z/m\mathbb Z)}
+O(m^3T^{2-\varepsilon}).
\end{aligned}
$$

We may replace $q_0$ by $6q_0$. Since $(m,6q_0)=1$ implies $(m,q_0)=1$, the congruence-counting estimate \eqref{eq:congruence-count} remains valid for every squarefree $m$ satisfying $(m,6q_0)=1$. Moreover, $(m,6q_0)=1$ implies $2\nmid m$ and $3\nmid m$, so every prime $\ell\mid m$ satisfies $\ell>3$. Hence the local density formula for $\beta(\ell)$ established earlier applies to every prime divisor $\ell$ of $m$. Relabeling $6q_0$ as $q_0$, we may assume $6\mid q_0$ throughout. So all sieve products below omit the prime divisors of this fixed $q_0$; by the finite-exception observations earlier, the sieve dimension and the power $-3$ in the local-factor product are unchanged.

Now for such a squarefree modulus $m$, define
$$
\beta(m)=\prod_{\ell\mid m}\beta(\ell),
\qquad \beta(1)=1,
$$
and let
$$
\mathcal B_m=
\{g\in\mathrm{SL}_2(\mathbb Z/m\mathbb Z):
F(g)\equiv0\pmod m\}.
$$

Since $m$ is squarefree, we can write $m=\ell_1\cdots\ell_k$, where $\ell_1,\ldots,\ell_k$ are distinct primes. By the Chinese remainder theorem,
$$
\mathbb Z/m\mathbb Z\cong\prod_{\ell\mid m}\mathbb F_\ell.
$$

For a matrix
$g=\begin{pmatrix}a&b\\c&d\end{pmatrix}$
with entries in $\mathbb Z/m\mathbb Z$, the condition $g\in\mathrm{SL}_2(\mathbb Z/m\mathbb Z)$ means that $ad-bc\equiv1\pmod m$. Since $m$ is squarefree, this is equivalent to $ad-bc\equiv1\pmod\ell$ for every prime $\ell\mid m$. Thus the reduction of $g$ modulo each $\ell\mid m$ belongs to $\mathrm{SL}_2(\mathbb F_\ell)$. Moreover, these reductions uniquely determine $g$ modulo $m$. Therefore,
$$
\mathrm{SL}_2(\mathbb Z/m\mathbb Z)
\cong\prod_{\ell\mid m}\mathrm{SL}_2(\mathbb F_\ell).
$$

Hence the condition $F(g)\equiv0\pmod m$ holds if and only if $F(g)\equiv0\pmod\ell$ for every prime $\ell\mid m$. By the definition of the local densities, we have
$$
\begin{aligned}
\frac{\#\mathcal B_m}
{\#\mathrm{SL}_2(\mathbb Z/m\mathbb Z)}
=
\prod_{\ell\mid m}
\frac{\#\{g\in\mathrm{SL}_2(\mathbb F_\ell):F(g)=0\}}
{\#\mathrm{SL}_2(\mathbb F_\ell)} =\prod_{\ell\mid m}\beta(\ell)
=\beta(m).
\end{aligned}
$$

Note that every matrix $\gamma\in\mathrm{SL}_2(\mathbb Z)$ satisfying $m\mid F(\gamma)$ belongs to exactly one residue class $g\in\mathcal B_m$ modulo $m$. Therefore, we can count these matrices by applying \eqref{eq:congruence-count} to each $g\in\mathcal B_m$ and summing the resulting counts. This gives
$$
\begin{aligned}
&\#\{\gamma\in\mathrm{SL}_2(\mathbb Z):
\|\gamma\|\leq T,\ m\mid F(\gamma)\} \\ 
&\qquad =
\frac{\#\mathcal B_m}
{\#\mathrm{SL}_2(\mathbb Z/m\mathbb Z)}M(T)
+O(\#\mathcal B_m\,m^3T^{2-\varepsilon})\\
&\qquad =
\beta(m)M(T)
+O(\#\mathcal B_m\,m^3T^{2-\varepsilon}).
\end{aligned}
$$
To bound the error, observe that
$$
\begin{aligned}
\#\mathcal B_m
\leq\#\mathrm{SL}_2(\mathbb Z/m\mathbb Z) =\prod_{\ell\mid m}\ell(\ell^2-1)
\leq m^3.
\end{aligned}
$$
Therefore,
\begin{equation}\label{eq:level-count}
\begin{aligned}
\#\{\gamma\in\mathrm{SL}_2(\mathbb Z):
\|\gamma\|\leq T,\ m\mid F(\gamma)\} =
\beta(m)M(T)+O(m^6T^{2-\varepsilon}).
\end{aligned}
\end{equation}

It remains to show that these errors are sufficiently small when summed over the sieve moduli. Let
$$
\begin{aligned}
e_m(T)
&=
\#\{\gamma\in\mathrm{SL}_2(\mathbb Z):
\|\gamma\|\leq T,\ m\mid F(\gamma)\} -\beta(m)M(T).
\end{aligned}
$$ be the error term.
By \eqref{eq:level-count}, we have
$$
|e_m(T)|\ll m^6T^{2-\varepsilon}.
$$

Hence, for any $L\geq1$,
$$
\begin{aligned}
\sum_{\substack{m\leq L\\
m\ {\rm squarefree}\\
(m,q_0)=1}}
|e_m(T)|
&\ll T^{2-\varepsilon}\sum_{m\leq L}m^6 \ll T^{2-\varepsilon}L^7.
\end{aligned}
$$

Now choose $L=T^\theta$, where $0<\theta<\min\{1,\varepsilon/7\}$, and let $\eta=\varepsilon-7\theta>0$. We have
$$
\begin{aligned}
\sum_{\substack{m\leq T^\theta\\
m\ {\rm squarefree}\\
(m,q_0)=1}}
|e_m(T)|
&\ll T^{2-\varepsilon+7\theta} =T^{2-\eta}.
\end{aligned}
$$

Note that $M(T)\asymp T^2$, as by taking $m=1$ in \eqref{eq:congruence-count}. Hence the total error $O(T^{2-\eta})$ is $o(M(T))$. Thus the accumulated congruence-counting error is smaller than the total number of matrices by a positive power of $T$, giving the required level of distribution for the sieve.

\medskip
\noindent\textit{(vii) Sieve setup.}
\par\nobreak\smallskip
We have now established the local densities, the sieve dimension $3$, and the required congruence distribution estimate. We next use these results to count matrices for which $F(\gamma)$ has no small prime divisors.

For $z\geq2$, let
$$
P(z)=
\prod_{\substack{3<\ell<z\\
\ell\ {\rm prime},\ \ell\nmid q_0}}\ell.
$$
We want to estimate the number of matrices $\gamma$ satisfying $\|\gamma\|\leq T$ and $(F(\gamma),P(z))=1$.

Before applying the sieve, we must exclude the matrices for which $F(\gamma)=0$, since the sieve sequence will consist of positive values of $|F(\gamma)|$. We first show that only finitely many such matrices exist.

Recall that
$$
F(\gamma)=R(\gamma)
\bigl(R(\gamma)-1-A(\gamma)\bigr)
\bigl(R(\gamma)-1+A(\gamma)\bigr).
$$
Since $R(\gamma)=x^2+3y^2>0$ for every $\gamma\in\mathrm{SL}_2(\mathbb Z)$, the condition $F(\gamma)=0$ implies
$$
A(\gamma)=\pm(R(\gamma)-1).
$$
From \eqref{eq:determinant-identity}, we have
$$
R(\gamma)B(\gamma)-A(\gamma)^2=3.
$$
Thus
$$
R(\gamma)\mid A(\gamma)^2+3
=(R(\gamma)-1)^2+3,
$$
which gives $R(\gamma)\mid4$.

Hence $R(\gamma)$ can take only finitely many positive values $1,2,4$. Since $R(\gamma)=x^2+3y^2$, each value permits only finitely many pairs $(x,y)$. For each such pair, the possible values $A(\gamma)=\pm(R(\gamma)-1)$ are also finite in number. Finally, the equations
$$
xu+3yv=A(\gamma),\qquad
xv-yu=1
$$
determine at most one pair $(u,v)$ for each fixed $(x,y,A(\gamma))$, because their determinant is $R(\gamma)>0$.

Therefore, only finitely many matrices satisfy $F(\gamma)=0$. More precisely, there is a constant $C_F>0$ such that $\#\{\gamma\in\mathrm{SL}_2(\mathbb Z):F(\gamma)=0\}\leq C_F$. In particular, the number of these matrices with $\|\gamma\|\leq T$ is at most $C_F$ for every $T$.

We can now define the nonnegative sequence to which we will apply the sieve. For $n\geq1$, let
$$
a_n(T)=
\#\{\gamma\in\mathrm{SL}_2(\mathbb Z):
\|\gamma\|\leq T,\ |F(\gamma)|=n\},
$$
and let
$$
\mathcal X(T)=\sum_{n\geq1}a_n(T).
$$
The sequence has finite support because $|F(\gamma)|\ll T^6$ whenever $\|\gamma\|\leq T$. By definition, $M(T)$ counts all matrices $\gamma\in\mathrm{SL}_2(\mathbb Z)$ with $\|\gamma\|\leq T$, whereas $\mathcal X(T)$ counts only those for which $F(\gamma)\neq0$. Thus,
$$
M(T)-\mathcal X(T)
=
\#\{\gamma\in\mathrm{SL}_2(\mathbb Z):
\|\gamma\|\leq T,\ F(\gamma)=0\}.
$$
That is, $\mathcal X(T)$ counts the matrices counted by $M(T)$ except those satisfying $F(\gamma)=0$.

We have shown that only finitely many matrices in $\mathrm{SL}_2(\mathbb Z)$ satisfy $F(\gamma)=0$. Hence the difference $M(T)-\mathcal X(T)$ is bounded independently of $T$, giving
$$
\mathcal X(T)=M(T)+O(1)\asymp T^2,
$$
where $M(T)\asymp T^2$ was shown by taking $m=1$ in \eqref{eq:congruence-count}.

\medskip
\noindent\textit{(viii) Verification of the sieve conditions.}
\par\nobreak\smallskip
We now verify the conditions needed to apply the upper-bound sieve to the sequence $a_n(T)$. We need to verify three conditions. First, for every squarefree $m$ with $(m,q_0)=1$, the number of matrices $\gamma$ satisfying $\|\gamma\|\leq T$, $F(\gamma)\neq0$, and $m\mid F(\gamma)$ must equal $\beta(m)\mathcal X(T)$ plus an error term $r_m(T)$ that satisfies a suitable bound when summed over the sieve moduli. Second, the function $\beta(m)$ must be multiplicative and give sieve dimension $3$. Third, for some $\theta,\delta>0$, the sum of $|r_m(T)|$ over all squarefree $m\leq T^\theta$ with $(m,q_0)=1$ must be $O(T^{2-\delta})$. We have already established the second condition, so it remains to verify the counting formula in the first condition and bound the total error in the third.

We first verify the counting formula. For every squarefree $m$ with $(m,q_0)=1$, we have already proved that
$$
\begin{aligned}
\#\{\gamma\in\mathrm{SL}_2(\mathbb Z):
\|\gamma\|\leq T,\ m\mid F(\gamma)\} =\beta(m)M(T)+O(m^6T^{2-\varepsilon}).
\end{aligned}
$$
The sum $\sum_{n\geq1,\,m\mid n}a_n(T)$ counts the same matrices, except those satisfying $F(\gamma)=0$. Since only finitely many matrices satisfy $F(\gamma)=0$, removing these matrices changes the count by at most $O(1)$. Hence
$$
\sum_{\substack{n\geq1\\m\mid n}}a_n(T)
=\beta(m)M(T)+O(m^6T^{2-\varepsilon}+1).
$$
Using $M(T)=\mathcal X(T)+O(1)$ and $0\leq\beta(m)\leq1$, we have
$$
\sum_{\substack{n\geq1\\m\mid n}}a_n(T)
=\beta(m)\mathcal X(T)+r_m(T),
$$
where $$|r_m(T)|\ll m^6T^{2-\varepsilon}+1.$$

Let $Z(T)=M(T)-\mathcal X(T)$. With $e_m(T)$ the remainder in the congruence count from (vi), the current remainder satisfies $r_m(T)=e_m(T)-(1-\beta(m))Z(T)$. Since $0\leq\beta(m)\leq1$ and $Z(T)=O(1)$ uniformly in $T$, we have $r_m(T)=e_m(T)+O(1)$ uniformly in $m$ and $T$.

It remains to verify the total error bound. We have already shown that the sum of the errors in counting matrices satisfying $m\mid F(\gamma)$ is $O(T^{2-\eta})$ over all squarefree $m\leq T^\theta$ with $(m,q_0)=1$, where $\eta=\varepsilon-7\theta>0$. Since excluding the finitely many matrices satisfying $F(\gamma)=0$ changes the error for each $m$ by at most $O(1)$, we have
$$
\begin{aligned}
\sum_{\substack{m\leq T^\theta\\
m\ {\rm squarefree}\\
(m,q_0)=1}}
|r_m(T)|
&\ll T^{2-\eta}+T^\theta \ll T^{2-\delta},
\end{aligned}
$$
where $\delta=\min\{\eta,2-\theta\}>0$, with $0<\theta<1$.

\medskip
\noindent\textit{(ix) Application of the upper-bound sieve.}
\par\nobreak\smallskip

We have now verified all three sieve conditions $(A_0)$--$(A_2)$ as formulated in \cite[Section~13.1]{BGS11}. Specifically, $(A_0)$ follows from the counting formula established above, the multiplicativity of $\beta$, and the bound $0\leq\beta(\ell)<1/2$ for all nonexceptional primes $\ell$; $(A_1)$ follows from the power-saving bound for the sum of the remainders at the level of distribution $D=L=T^\theta$; and $(A_2)$ follows from the dimension-three estimate after excluding the finitely many exceptional primes. Since $\mathcal X(T)\asymp T^2$ and $\theta<1$, we also have $D<\mathcal X(T)$ for all sufficiently large $T$. We can therefore apply the upper-bound sieve~\cite[Theorem~13.1]{BGS11} with sieve dimension $3$.

Choose any $s>27$ and let
$$
z=L^{1/s}=T^\nu,
\qquad
\nu=\frac{\theta}{s}.
$$
Since $0<\theta<1$ and $s>27$, we have $0<\nu<1/2$.

The upper-bound sieve gives
$$
\sum_{\substack{n\geq1\\(n,P(z))=1}}a_n(T)
\ll
\frac{\mathcal X(T)}
{(\log\mathcal X(T))^3}
\ll\frac{T^2}{(\log T)^3}.
$$
It follows that
\begin{equation}\label{eq:matrix-sieve-bound}
\#\{\gamma\in\mathrm{SL}_2(\mathbb Z):
\|\gamma\|\leq T,\ (F(\gamma),P(T^\nu))=1\}
\ll\frac{T^2}{(\log T)^3}.
\end{equation}

Note that for sufficiently large $T$, the integer $P(T^\nu)$ is greater than $1$. Hence matrices with $F(\gamma)=0$ do not satisfy $(F(\gamma),P(T^\nu))=1$, so the set on the left agrees with the sifted sequence above.

\medskip
\noindent\textit{(x) Counting trey primes.}
\par\nobreak\smallskip
We now return to the original problem of counting trey primes. To use \eqref{eq:matrix-sieve-bound}, we need to show that each trey prime in the range considered below gives a distinct matrix counted by the left-hand side.

We first separate the trey primes $r\leq\sqrt X$. There are at most $\sqrt X$ such primes, so $N_3(\sqrt X)\leq\sqrt X.$ We will count the remaining trey primes $\sqrt X<r\leq X$.

For each such prime, choose the matrix $\gamma$ constructed above from its trey group $\{1,p,q\}$. Since $\|\gamma\|^2<2r\leq2X$, setting $T=\sqrt{2X}$ gives $\|\gamma\|<T$.

We next show that $F(\gamma)$ is coprime to $P(T^\nu)$. By \eqref{eq:F-trey}, we have $F(\gamma)=4rpq.$ Thus it is enough to show that the three primes $r$, $p$, and $q$ are all greater than $T^\nu$, since $P(T^\nu)$ contains neither $2$ nor any prime dividing $q_0$.

We already have $r>\sqrt X>X^{1/4}$. To obtain the corresponding bounds for $p$ and $q$, recall that the trey group relations give
$$
p^2\equiv q\pmod r,\qquad
q^2\equiv p\pmod r.
$$
If $p<\sqrt r$, then $p^2<r$. Since $q$ is the standard representative of $p^2$ modulo $r$, this would imply $q=p^2$, contradicting the primality of $q$. Hence $p>\sqrt r$. The same argument, with $p$ and $q$ interchanged, gives $q>\sqrt r$.

Therefore,
$$
p,q>\sqrt r>X^{1/4}.
$$
Since $\nu<1/2$ and $T=\sqrt{2X}$, we have $T^\nu=(2X)^{\nu/2}<X^{1/4}$ for sufficiently large $X$.

Thus every prime factor of $F(\gamma)=4rpq$ other than $2$ exceeds $T^\nu$. It follows that
$$
(F(\gamma),P(T^\nu))=1.
$$
Hence every trey prime $\sqrt X<r\leq X$ gives a matrix counted by \eqref{eq:matrix-sieve-bound}.

Moreover, distinct trey primes give distinct matrices, since the value $R(\gamma)=r$ determines the prime modulus. Therefore,
$$
\begin{aligned}
N_3(X)-N_3(\sqrt X)
&\leq
\#\{\gamma\in\mathrm{SL}_2(\mathbb Z):
\|\gamma\|\leq T,\ (F(\gamma),P(T^\nu))=1\} \ll\frac{T^2}{(\log T)^3}.
\end{aligned}
$$
Since $T=\sqrt{2X}$, we have
$$
N_3(X)-N_3(\sqrt X)
\ll\frac{X}{(\log X)^3}.
$$

Finally,
$$
N_3(\sqrt X)\leq\sqrt X
\ll\frac{X}{(\log X)^3}.
$$

Combining the two bounds gives
$$
N_3(X)\ll\frac{X}{(\log X)^3}
$$
for all sufficiently large $X$, say $X\geq X_0$. To complete the claim for $X\geq2$, we can always enlarge the implied constant so that the same estimate also holds for the remaining range $2\leq X<X_0$. This completes the proof.

\end{proof}

\begin{remark}
The relation $r\mid(q-p)^2+3$ also gives a restriction on the difference between the two prime elements of a trey group. Writing $a=q-p>0$, we have $a^2+3\geq r$, so $a\geq\sqrt{r-3}$. Moreover, for any fixed positive integer $a$, the prime modulus $r$ must divide $a^2+3$, which has only finitely many prime divisors. Since each prime modulus has at most one subgroup of order $3$, only finitely many trey groups can have a fixed gap $a$.
\end{remark}

We finish this section by associating the number of trey primes with the number of primes congruent to $1$ modulo $3$.

\begin{corollary}\label{cor:relative-density}
Among primes $r\equiv1\pmod3$ with $r\leq X$, the proportion for which both standard roots of $x^2+x+1\equiv0\pmod r$ are prime is $O((\log X)^{-2})$. In particular, this proportion tends to zero as $X\to\infty$.
\end{corollary}

\begin{proof}
By Corollary~\ref{cor:prime-roots}, the primes $r\equiv1\pmod3$ for which both standard roots are prime are exactly the trey primes. Thus their number up to $X$ is $N_3(X)$.

By the prime number theorem in arithmetic progressions \cite{Davenport2000},
$$
\pi(X;3,1)\sim\frac{X}{2\log X}.
$$
Combining this with Theorem~\ref{thm:main-upper}, we obtain
$$
\frac{N_3(X)}{\pi(X;3,1)}
\ll
\frac{X/(\log X)^3}{X/\log X}
=\frac1{(\log X)^2}.
$$
Since $(\log X)^{-2}\to0$ as $X\to\infty$, the required proportion tends to zero.
\end{proof}

\section{Open problems}\label{sec:questions}

We have proved that $N_3(X)\ll\frac{X}{(\log X)^3}$, but this upper bound does not determine how frequently trey primes occur. A natural question is whether $N_3(X)$ has the same order of magnitude as this upper bound.

\begin{question}\label{q:asymptotic}
Is there a constant $C_3>0$ such that $N_3(X)\sim C_3\frac{X}{(\log X)^3}$?
\end{question}

Such an asymptotic formula would imply that there are infinitely many trey primes. However, even the existence of infinitely many trey primes remains unproved without Schinzel's Hypothesis H.

\begin{question}\label{q:unconditional}
Can one prove unconditionally that $N_3(X)\longrightarrow\infty$, or prove that $N_3(X)$ is bounded which would disprove Schinzel's Hypothesis H?
\end{question}

We can also ask whether groups of prime units of larger orders occur infinitely often. For a prime modulus $r$, a subgroup of order $d$ can exist only if $d\mid r-1$. The nonidentity elements of such a subgroup are roots of $\frac{x^d-1}{x-1}$ modulo $r$, and the elements of exact order $d$ are roots of the $d$-th cyclotomic polynomial $\Phi_d(x)$. For the subgroup to be a group of prime units, all its nonidentity elements must have prime standard representatives. This leads to the following question.

\begin{question}\label{q:fixed-order}
For which odd integers $d\geq3$ does $N_d(X)\longrightarrow\infty$?
\end{question}

A related but different question is whether the orders of groups of prime units modulo primes can be arbitrarily large. This does not require any particular fixed order to occur infinitely often. Since every subgroup of a group of prime units is itself a group of prime units, Cauchy's theorem implies that a group whose order is divisible by a prime $\ell$ contains a group of prime units of order $\ell$.

\begin{question}\label{q:unbounded-order}
Do groups of prime units modulo primes have unbounded order?
\end{question}

We then consider groups of prime units modulo general integers. We have proved that, under Dickson's conjecture, every fixed order $d\geq2$ occurs for infinitely many composite moduli. In particular, the orders of these groups are unbounded under Dickson's conjecture. We ask whether the same conclusion can be proved unconditionally.

\begin{question}\label{q:general-unbounded-order}
Can one prove unconditionally that groups of prime units modulo general integers have unbounded order, or prove that their orders are bounded which would disprove Dickson's conjecture?\end{question}

The above questions concern groups of prime units modulo integers. It is also natural to ask whether analogous results hold over polynomial rings over finite fields. To formulate this question, one possible convention is to work modulo a nonconstant polynomial $f\in\mathbb F_q[t]$, choose the unique representative of degree less than $\deg f$ in each residue class, and call a nonconstant representative prime when it is irreducible over $\mathbb F_q$.

With this convention in mind, we end this paper with a final question.

\begin{question}
To what extent do the results and open questions of this paper extend to $\mathbb F_q[t]$, with prime numbers replaced by irreducible polynomials?
\end{question}

\end{document}